\documentclass[11pt,a4paper]{article}
\usepackage{amsmath,amssymb,amsthm,mathrsfs,fullpage,lineno,bbm}
\usepackage{hyperref}

\newtheorem{theorem}{Theorem}
\newtheorem{lemma}{Lemma}

\theoremstyle{definition}
\newtheorem{definition}{Definition}

\newtheorem{remark}{Remark}
\newtheorem{example}{Example}

\newtheorem{fact}{Fact}

\newtheorem{ques}{Question}
\newcommand{\E}{\mathbb{E}}

\newcommand\tup[1]{\left\langle #1 \right\rangle}

\title{Functional inequalities and random walks on increasing subsets of the hypercube}

\begin{document}
\author{Fan Chang\thanks{School of Statistics and Data Science, Nankai University, Tianjin, China and Extremal Combinatorics and Probability Group (ECOPRO), Institute for Basic Science (IBS), Daejeon, South Korea. E-mail: 1120230060@mail.nankai.edu.cn.  Supported by the NSFC under grant 124B2019 and the Institute for Basic Science (IBS-R029-C4).}, Guowei Sun\thanks{School of Mathematics, Shandong University, Jinan, China. E-mail: gwsun@mail.sdu.edu.cn} and Lei Yu\thanks{School of Statistics and Data Science, LPMC, KLMDASR, and LEBPS, Nankai University, Tianjin, China. E-mails: leiyu@nankai.edu.cn. Supported by the National Key Research and Development Program of China under grant 2023YFA1009604, the NSFC under grant 62101286, and the Fundamental Research Funds for the Central Universities of China (Nankai University) under grant 054-63253112.}}
\date{}

\maketitle

\begin{abstract}
Motivated by random walks on subsets of the hypercube, we prove two discrete functional inequalities on the hypercube by the technique of induction-by-restrictions. First, we give a short, elementary proof of the Poincar\'e inequality on increasing subsets of the cube recently established by Fei and Ferreira Pinto Jr~\cite{fei2025spectral}, which yields an $O(n^2)$ upper bound on the mixing time of censored random walks, improving upon previous bounds. Second, adapting Samorodnitsky's induction method~\cite{Samorodnitsky17} to the $p$-biased setting, we establish a sharp $p$-biased edge-isoperimetric inequality for real-valued functions supported on increasing sets, which recovers the classic biased edge-isoperimetric inequality for increasing sets and identifies increasing subcubes as the extremizers. This result also admits a probabilistic interpretation in terms of maximizing the mean first exit time of biased random walks known as Glauber dynamics. 

\vspace{10pt}

\vspace{10pt}

\end{abstract}

\section{Introduction}
Functional inequalities on the discrete hypercube/Hamming graph, such as the Poincar\'e, log-Sobolev, and edge-isoperimetric inequalities, form a central part of modern discrete analysis. They connect Boolean function analysis, discrete isoperimetric problems, and the spectral theory of Markov chains, providing powerful tools for studying concentration of measure, threshold phenomena, and mixing times of random walks. In the classical uniform product setting on the hypercube $\{0,1\}^n$, these inequalities are well understood: sharp constants are known (e.g., the exact log-Sobolev constant~\cite{BT2006logsobolev}) and elegant proofs arise via tensorization~\cite{gross1975logarithmic}, semigroup methods~\cite{Ivanisvili2024KKL,ivanisvili2020rademacher}, or discrete Fourier analysis~\cite{Eldan2025Isoperimetric}. In this paper, we depart from the uniform product setting---by restricting to structured subsets, specifically increasing subsets, or by working under biased measures. We aim at deriving sharp functional inequalities in these spaces and applying them to characterize the convergence behaviors  of  relevant random walks.
\subsection{Poincar\'e inequality and random walks on increasing sets}
The study of random walks censored on increasing subsets of the hypercube is motivated by foundational challenges in Markov Chain Monte Carlo (MCMC) for constrained sampling~\cite{guruswami2016rapidly,Jer1998,LP2017,MT2006}. A canonical example is sampling uniformly from configurations of a statistical mechanics model (e.g., critical percolation~\cite{DM2014}) that satisfy a non-local property such as a left-to-right crossing. While the ambient space is equipped with a product measure and can be easily sampled, the associated constraint set forms a highly correlated, nonproduct subset. Rejection sampling becomes statistically inefficient for rare events, prompting the use of local chains that \emph{remain within} the constraint set. The natural censored dynamics re-randomizes one coordinate at a time and accepts the move only if the property is preserved. When the property is monotone, the state space forms an increasing (upward-closed) subset, where this monotonicity ensures the irreducibility of the chain. The central algorithmic question is to quantify the mixing time and the related survival/exit characteristics of this censored walk, since these quantities determine the efficiency and robustness of constrained sampling in combinatorics, statistical physics, and property testing.

We work on the $n$-dimensional hypercube, more precisely,  the Hamming graph. The Hamming graph is the graph $Q^n=(\{0,1\}^n,E)$ with $2^n$ vertices indexed by binary strings, where two vertices are adjacent (written $x\sim y$) iff they differ in exactly one coordinate.  Discrete-time random walks on the Hamming graph  (or more generally, on  finite graphs $G=(V,E)$)  model stochastic processes where a walker iteratively transitions between adjacent vertices. Formally, at each timestep $t$, the walker moves from the current vertex $u$ to a neighbor $v$ with probability $P(u,v)=\frac{1}{{\deg}(u)}$, defining a Markov chain with state space $V$ and transition matrix $P$. The stationary distribution $\pi$ of this random walk, satisfying $\pi P=\pi$, is given by $\pi(v)=\frac{\deg(v)}{2|E|}$ for connected non-bipartite graphs. The mixing time $t_{{\rm mix}}(\varepsilon)$ quantifies the convergence rate to stationarity, defined as the minimal $t$ such that:
\begin{equation}
    \max_{u\in V}\|P^{t}(u,\cdot)-\pi\|_{{\rm TV}}\le \varepsilon,
\end{equation}
where $\|\cdot\|_{{\rm TV}}$ denotes the total variation, and $\varepsilon\in(0,1)$ governs approximation accuracy. 
 
Ding and Mossel~\cite{DM2014} initiated the study of the \emph{random walk on $\{0,1\}^n$ censored to a subset $A \subset \{0,1\}^n$}. From any state $x \in A$, the walk attempts to flip a coordinate chosen uniformly at random; it transitions to the resulting neighbor if it remains within $A$, otherwise it stays at $x$. Notably, without imposing specific restrictions on $A$, the censored random walk can exhibit widely varying behavior. Even when $A$ is both large and connected, the walk may still mix extremely poorly (see [14, Examples 1 and 2]). A significant open question (\cite[Question 1.1]{DM2014}) asked if the optimal mixing time belongs to $O_c(n\log n)$ for increasing subsets $A$ with constant density under the uniform measure. Here, a subset $A\subseteq\{0,1\}^n$ is called \emph{increasing} if for all $x,y\in\{0,1\}^n$,
\begin{equation*}
x_i\le y_i,\ \forall i\in[n], \quad x\in A \ \Rightarrow\ y\in A.
\end{equation*}

Let $\mu$ be the uniform probability measure on $\{0,1\}^n$. 
\begin{ques}[{\cite[Question 1.1]{DM2014}}]
	Let $\varepsilon,c \in (0,1)$ be  constants. For an increasing subset $A\subseteq\{0,1\}^n$ with measure $\mu(A)\ge c$, consider the censored random walk on $A$ defined above. Is it true that $t_{{\rm mix}}(\varepsilon)\le O_{\varepsilon,c}(n\log n)$?     
\end{ques}

When the density $\mu(A)$ is constant, Ding and Mossel established a mixing time bound of $O(n^3)$. Very recently, Fei and Ferreira Pinto Jr.~\cite{fei2025spectral} made a breakthrough on this topic by showing an improved mixing time bound of $O(n^2)$ for constant-density increasing sets. Their result is obtained via establishing the following Poincar\'e inequality on the non-product subspace. We write ${\rm Var}_A[f]$ for the variance of $f(x)$ when $x$ is drawn uniformly from $A$. Moreover, $\mathcal{E}_A(f)$ denotes the ``restricted'' Dirichlet form of $f$ on $A$, defined by
\begin{equation}\label{eq:restricted-Dirichlet-form}
\mathcal{E}_A(f)=\mathcal{E}_A(f,f):=\frac{1}{4}\sum_{i=1}^n \underset{x\sim A}{\mathbb{E}}\!\left[(f(x)-f(x^{\oplus i}))^2\,\mathbbm{1}_{\{x^{\oplus i}\in A\}}\right],
\end{equation}
where the expectation is taken with respect to the uniform probability measure on the increasing subset $A\subseteq \{0,1\}^n$ and $x^{\oplus i}$ is obtained from $x$ by flipping the $i$-th coordinate.
\begin{theorem}[{\cite[Theorem 1.6]{fei2025spectral}}]\label{thm:main F25}
Let $A\subseteq \{0,1\}^n$ be a non-empty increasing set. It holds that for all $f:A\to \mathbb{R}$,
\begin{equation}
    {\rm Var}_A[f]\le\frac{1}{1-\sqrt{1-\mu(A)}}\cdot \mathcal{E}_A(f).
\end{equation}
\end{theorem}

When $A=\{0,1\}^n$, the above result recovers the classical Poincar\'e inequality on the hypercube. The proof of this theorem provided by Fei and Ferreira Pinto Jr relies heavily on the connections between the directed heat flow and the directed Poincar\'e inequality. In this paper, we provide a simple inductive proof of their Poincar\'e inequality, though yielding a slightly worse constant (theirs being 1 versus ours at 2). Specifically, we use the induction method to prove the following theorem, a weaker version of Theorem~\ref{thm:main F25}. 
\begin{theorem}\label{thm:main thm2}
Let $A\subseteq \{0,1\}^n$ be a non-empty increasing set. Then for all $f:A\to \mathbb{R}$,
\begin{equation}
    {\rm Var}_A[f]\le\frac{2}{1-\sqrt{1-\mu(A)}}\cdot \mathcal{E}_A(f).
\end{equation}
\end{theorem}

Nevertheless, we still obtain an $O(n^2)$ mixing time bound for increasing sets $A$ of constant density via standard Markov chain theory (\cite[Theorem 12.4]{LP2017}). 
\begin{theorem}
Let $A\subseteq\{0,1\}^n$ be a non-empty increasing set. Then the random walk on $\{0,1\}^n$ censored to $A$ has mixing time
$$
t_{{\rm mix}}\le \frac{2n}{\mu(A)}\cdot \log(4\cdot 2^n\mu(A)).
$$
\end{theorem}

Fei and Ferreira Pinto Jr~\cite{fei2025spectral} develop an $L^2$ functional–analytic framework for increasing sets $A\subseteq\{0,1\}^n$. They introduce a \emph{directed} gradient that counts only order-consistent edges (down moves), which defines a directed Laplacian $L_A^{\to}$ and energy $\mathcal E_A^{\to}(f)$. They then run the monotone heat flow generated by $L_A^{\to}$ and prove an energy–variance dissipation inequality along the flow, yielding a \emph{directed Poincar\'e inequality} with constant $C(\mu(A))=O(1)$ whenever $\mu(A)$ is bounded away from $0$. A second ingredient is an \emph{approximate FKG inequality} on increasing sets that supplies quantitative positive-correlation bounds. Together these tools give an optimal (up to constants) spectral gap for the censored walk on $A$. 

In contrast, our method is designed to be simple and straightforward. We establish a slightly weaker, yet asymptotically optimal, Poincar\'e inequality by a classical induction-by-restrictions argument. The proof restricts functions to lower-dimensional subcubes and uses the inherent increasing structure to decompose the variance and the Dirichlet form fiberwise, ultimately reducing the task to an explicit five-point inequality that can be verified directly. The resulting constant is within a factor $2$ of their value and implies the same mixing bound $O(n^2)$ for censored walks. Because every step is combinatorial and closed-form, the proof is short, self-contained, and readily adaptable.

\subsection{Biased Samorodnitsky's inequality and exit time of random walks from increasing sets}

Another motivation for studying functional inequalities supported on non-product sets is to characterize the mean first exit time of a random walk from an increasing set.

For functions $g,h:\{0,1\}^n\to\mathbb{R}$, we define the \emph{Dirichlet form} as
$$
\mathcal{E}(g,h)=\frac{1}{4}\underset{x\sim\mu}{\mathbb{E}}\left[\sum\limits_{y\in\{0,1\}^n:y\sim x}(g(x)-g(y))(h(x)-h(y))\right],
$$
and write $\mathcal{E}(f):=\mathcal{E}(f,f)$, where the expectation is taken with respect to the uniform probability measure $\mu:=\mu_{1/2}^{n}$ on the hypercube. In 2017, Samorodnitsky~\cite{Samorodnitsky17} established the following discrete functional inequality on the hypercube.
\begin{theorem}[Samorodnitsky~\cite{Samorodnitsky17}]\label{thm:Samorodnitsky main thm}
Let $A$ be a non-empty subset of $\{0,1\}^n$ and let $g:\{0,1\}^n\to \mathbb{R}$ be supported on $A$. Then
\begin{equation}\label{ineq: Samorodnitsky functional edge-isoper}
\mathcal{E}(g)\ge \frac{|A|}{2^{n+1}}\left(\underset{x\sim A}{\mathbb{E}}[|g(x)|]\right)^2\cdot\log_2\left(\frac{2^n}{|A|}\right),
\end{equation}
with equality if and only if $g$ is the indicator function of a co-dimension $k$ subcube for any $0\le k\le n$.
\end{theorem}

When $g$ is the indicator function of $A$, Samorodnitsky's inequality~\eqref{ineq: Samorodnitsky functional edge-isoper} implies the classic edge-isoperimetric inequality on $Q^n$ (see~\cite{Bernstein67,Harper64,Hart76,Lindsey64}): for every $A\subset\{0,1\}^n$,
\begin{equation}\label{ineq: Harper origin edge-isop}
|\nabla A|\ge |A|\cdot\log_2\left(\frac{2^n}{|A|}\right),
\end{equation}
where the \emph{edge boundary} $\nabla(A)$ is the set of edges between $A$ and its complement $A^c=\{0,1\}^n\setminus A$. Equality in~\eqref{ineq: Harper origin edge-isop} holds if and only if $A$ is a subcube of co-dimension $k$, for any $0\le k\le n$. 

Besides the motivation from generalization of  Harper's edge-isoperimetric inequality,  Samorodnitsky's inequality admits  a natural interpretation in the language of random walks on the hypercube. Specifically, it is equivalent to a sharp estimate of mean first exit time of a random walk from a subset with a given size; see~\cite{Samorodnitsky17} for further details. 

We now generalize Samorodnitsky’s inequality to the $p$-biased setting so that the resulting functional inequality recovers the sharp version of the $p$-biased edge-isoperimetric inequality. Let $f\in L^2(\{0,1\}^n,\mu_p^{n})$, where $\mu_p$ is the Bernoulli measure on $\{0,1\}$ with $\mu_p(1)=p$ and $\mu_p(0)=1-p$. When there is no ambiguity, we omit the superscript and use $\mu_p$ to denote $\mu_p^{n}$. For $x\in \{0,1\}^n$, write $|x|=|\{i\in[n]:x_i=1\}|$, so that
\begin{equation}\label{def:pbiased}
\mu_p(x):=\mu_p^{n}(x)=p^{|x|}(1-p)^{n-|x|}, \quad \mathbb{E}_p[f]:=\underset{x\sim\mu_p}{\mathbb{E}}[f(x)]=\sum\limits_{x\in \{0,1\}^n}f(x)\mu_p(x).
\end{equation}
Studying biased measures on the discrete cube is natural, since many applications (e.g., to the study of percolation~\cite{BKS1999}, threshold phenomena in random graphs~\cite{KLMM2024globalhyper2024,Talagrand1994randomgraph}, and hardness of approximation~\cite{DS2005Annal}) rely on their use. 

The definition of Dirichlet quadratic form of $g$ with respect to the biased measure is
\begin{equation}\label{eq:biased Dirichlet}
\mathcal{E}_p(g)=\mathcal{E}_p(g,g):=\frac{1}{4}\cdot\underset{x\sim\mu_p}{\mathbb{E}}\left[\sum\limits_{y\in\{0,1\}^n:y\sim x}(g(x)-g(y))^2\right].
\end{equation}
With this convention, each undirected edge is counted once (up to the normalization), which is consistent with the restricted Dirichlet form $\mathcal E_A$ used in Theorem~\ref{thm:main F25}.

The following `biased' generalization of Harper's edge-isoperimetric inequality for increasing subsets is considered folklore (see~\cite{EKL2019biased,Kahnkalai07}).
\begin{theorem}
Let $A$ be an increasing subset of $\{0,1\}^n$ and let $0 < p < 1$. Then
\begin{equation}\label{ineq: p-biased isoper}
4p\cdot\mathcal{E}_p(\mathbbm{1}_A)\ge \mu_p(A)\log_p\mu_p(A),
\end{equation}
with equality if and only if $A$ is an increasing subcube. 
\end{theorem}

Inspired by \eqref{ineq: Samorodnitsky functional edge-isoper}, it is natural to seek a functional inequality on the $p$-biased cube for general real-valued (increasing) functions such that the resulting functional inequality, in the special case of indicator functions, reduces to the $p$-biased edge-isoperimetric inequality with the \emph{sharp} constant. The standard log-Sobolev inequality is too weak to satisfy this requirement. We now establish a desired inequality analogous to Samorodnitsky's inequality in Theorem~\ref{thm:Samorodnitsky main thm}, which recovers \eqref{ineq: p-biased isoper} by substituting $g=\mathbbm{1}_A$.

\begin{theorem}\label{thm:main thm1}
Let $A\subseteq \{0,1\}^n$ be a non-empty increasing set, let $0<p<1$, and let $g:\{0,1\}^n\to\mathbb{R}$ be supported on $A$. Then
\begin{equation}\label{ineq: Samorodnitsky-type p-biased isoper}
4p\cdot\mathcal{E}_p(g)\ge \frac{\mathbb{E}_p[|g|]^2}{\mu_p(A)}\log_p\mu_p(A).
\end{equation}
Moreover, equality holds in \eqref{ineq: Samorodnitsky-type p-biased isoper} if and only if $A$ is an increasing subcube and $g$ is a scalar multiple of $\mathbbm{1}_A$.
\end{theorem}

We view the $p$-biased hypercube as a weighted undirected graph $Q_p^n$, where for every vertex $x\in\{0,1\}^n$ and coordinate index $i\in[n]$, the edge weight is defined by
\begin{equation*}
w_p(\{x,x^{\oplus i}\})=\mu_p(x)+\mu_p(x^{\oplus i})
=\mu_p(x_{\backslash i}),
\end{equation*}
where $x_{\backslash i}=(x_1,\dots,x_{i-1},x_{i+1},\dots,x_n)$. 
Thus each edge is weighted proportionally to the total $p$-biased mass of its two endpoints.

A random walk on a graph is a process that begins at some vertex and, at each time step, moves to another vertex according to a prescribed transition rule. When the graph is unweighted (i.e., $p=1/2$), the vertex visited by the (simple) random walk is chosen uniformly at random among the neighbors of the current vertex. This is just the random walk considered by Samorodnitsky~\cite{Samorodnitsky17}. However, any extension of this random walk to the $p$-biased setting seems unnatural. In the $p$-biased setting, there are several natural choices of dynamics in the literature. For our purposes, the most natural one is the \emph{Glauber} (or \emph{heat-bath}) dynamics, since the associated Dirichlet form is  exactly $\mathcal{E}_p(g)$ given in~\eqref{eq:biased Dirichlet}. 

More precisely, from a state $x\in\{0,1\}^n$, choose a coordinate $i\in[n]$ uniformly at random and resample the $i$-th coordinate according to the law $\mu_p$. Equivalently, the transition probability from $x$ to $x^{\oplus i}$ is given by 
\[
\mathbf P_p^{\mathrm{Gl}}(x,x^{\oplus i})
=\frac{\mu_p(1-x_i)}{n},
\quad i\in[n],
\]
and with the remaining probability, the walk stays at $x$. If the walk starts from a random point drawn from the product measure $\mu_p$, then the $i$-th coordinate is reset from $0$ to $1$ with probability $p$, while it is reset from $1$ to $0$ with probability $1-p$. In particular, this Markov chain is reversible, and hence stationary, with respect to $\mu_p$.

Now let $A\subseteq\{0,1\}^n$ with $\mu_p(A)>0$, and choose $X_0\in A$ according to the conditional law $\nu=\mu_p(\cdot|A)$, given by  $\nu(x):=\frac{\mu_p(x)}{\mu_p(A)},\forall x\in A$. Let $(X_t)_{t\ge 0}$ be the Glauber dynamics starting from $X_0\sim \nu$, and let
\[
\tau:=\min\{t\ge 1:X_t\notin A\}
\]
be the first time that the walk exits $A$. We refer to $\mathbb E_\nu[\tau]$ as the \emph{mean first exit time} of $A$ under the $p$-biased Glauber dynamics.

In~\cite{Samorodnitsky17}, Samorodnitsky observed that the functional inequality in~\eqref{ineq: Samorodnitsky functional edge-isoper} is equivalent to a probabilistic statement about random walks on the hypercube: among all subsets of a given cardinality, subcubes maximize the \emph{mean first exit time}. Our new $p$-biased edge-isoperimetric inequality in~\eqref{ineq: Samorodnitsky-type p-biased isoper} admits an analogous probabilistic interpretation in terms of the Glauber dynamics on the $p$-biased hypercube.
\begin{theorem}[Mean first exit time]\label{thm:mean-first-exit-time}
Let $A\subsetneq\{0,1\}^n$ be a non-empty increasing set and let $0<p<1$.
Let $(X_t)$ be the $p$-biased Glauber dynamics starting from $X_0\sim\nu=\mu_p(\cdot\mid A)$, and
let $\tau$ be the first exit time from $A$. Then
\[
\E_\nu[\tau] \le \frac{n}{(1-p)\log_p\bigl(\mu_p(A)\bigr)}.
\]
Moreover, equality holds if and only if $A$ is an increasing subcube.
\end{theorem}

\begin{example}[Increasing subcubes]\label{exam:subcube-exit}
Let $A=\{x\in\{0,1\}^n:x_1=\cdots=x_k=1\}$. Conditional on the chain being in $A$ at the current step, it exits at the next step iff the chosen coordinate lies in $[k]$ and it is resampled to $0$, which has probability $\frac{k}{n}(1-p)$. Hence $\tau$ is geometric on $\{1,2,\dots\}$ with parameter $\frac{k}{n}(1-p)$ and
\[
\E_\nu[\tau]=\frac{n}{k(1-p)}.
\]
Since $\mu_p(A)=p^k$, we have $\log_p(\mu_p(A))=k$, and thus, Theorem~\ref{thm:mean-first-exit-time} is tight
for subcubes.
\end{example}

Note that the random walk we consider differs from that in Samorodnitsky’s work: in his setting, the walk is the simple random walk on the cube (without staying at the current state $x$). In contrast, when specialized to the unbiased case $p=1/2$, our Glauber dynamics corresponds to the lazy simple random walk, where the walk stays at the current state $x$ with probability $1/2$. Nevertheless, our Theorem~\ref{thm:mean-first-exit-time} is consistent with his in the following sense. In our setting, the time to move from a state to an adjacent state is a geometric random variable with mean $2$. By contrast, this transition time is identically $1$ in Samorodnitsky’s setting. Accordingly, our estimate in Theorem~\ref{thm:mean-first-exit-time} with $p=1/2$ is exactly double the one given in~\cite[Theorem 1.2]{Samorodnitsky17}.

In the language of statistical mechanics, Theorem~\ref{thm:mean-first-exit-time} may be viewed as a result describing the temporal stability of a macroscopic state under stochastic thermal fluctuations in a system of non-interacting components. The configuration space is represented by the $n$-dimensional hypercube $\{0,1\}^n$, and the Glauber dynamics describes the system's evolution toward equilibrium through local spin-flip processes driven by an external field (represented by the bias $p$). When we define an increasing set $A$ as a ``macrostate'' (such as a ``connected cluster'' or an ``ordered phase''), the mean first exit time $\mathbb{E}_\nu[\tau]$ can be viewed as the metastable lifetime of that state. A pivotal finding in Theorem~\ref{thm:mean-first-exit-time} is that increasing subcubes are the unique extremizers that maximize the exit time. Physically, this implies that order maintained by fixed, local constraints (where specific spins are ``locked'' into a state) is the most resilient against thermal noise in the absence of internal coupling. Compared to any other macrostate of the same statistical weight ($\mu_p(A)$), the subcube structure minimizes its ``leakage'' into the disordered complement, thereby maximizing its survival time. Our upper bound on $\mathbb{E}_\nu[\tau]$ quantifies how the environment’s bias ($p$) competes with the system's geometric structure to determine how long an ordered state can persist before a sufficiently large fluctuation triggers a transition to the disordered phase.

\subsection{Organization}
This paper is organized as follows. In section~\ref{sec:monotone}, we present an inductive proof of the Poincar\'e inequality on increasing subsets of the cube given in Theorem~\ref{thm:main thm2}. Our arguments follow the induction-plus-five-point inequality paradigm, with the associated analytic verification being considerably simpler. Section~\ref{sec:biased} presents an inductive proof of the sharp $p$-biased edge-isoperimetric inequality for real-valued functions supported on increasing sets given in Theorem~\ref{thm:main thm1}. We develop an induction-by-restrictions framework and verify the required two-point inequality within this setting. Section~\ref{sec:mean first exit time} demonstrates how inequality~\eqref{ineq: Samorodnitsky-type p-biased isoper} in Theorem~\ref{thm:main thm1} yields bounds on the mean survival (first-exit) time of the $p$-biased Glauber dynamics on increasing sets (Theorem~\ref{thm:mean-first-exit-time}). Finally, section~\ref{sec:discuss} outlines a direction for future work: establishing suitable (possibly modified) log-Sobolev inequalities on increasing sets, which are closely related to the Ding--Mossel conjecture.

\section{Proof of Theorem~\ref{thm:main thm2}}\label{sec:monotone}
\subsection{Preliminaries}\label{sec:Prethm2}
In this subsection, we outline some notation and terminology, and present some basic lemmas for inductive proof of Theorem~\ref{thm:main thm2}. We use the so-called induction-by-restrictions method, a standard technique in combinatorics and discrete analysis (see, e.g.,~\cite{alpay2025lowerboundsdyadicsquare,beltran2023sharpisoperimetricinequalitieshypercube,BIMP2025young,CTKS2025inequalities,durcik2024sharpisoperimetricinequalitieshamming,EKL2019biased,FS2007CPCKKL,Kahnkalai07,Samorodnitsky17,yu2025averagedistancelevel1fourier}), and proceeds by considering the restricted functions $f_z$ obtained by fixing the last coordinate to $z\in\{0,1\}$. One then expresses Dirichlet forms (or total influence) and variances recursively in terms of $f_0$ and $f_1$, applies the inductive hypothesis in dimension $n-1$, and controls the remaining cross-terms by verifying an explicit finite-point inequality. In fact, there is another standard inductive method in the literature which proceeds  induction-by-derivatives; this method underlies classic proofs of Bonami's lemma and the hypercontractive inequality~\cite{Ryanbook2014}.

Given $A\subset\{0,1\}^{n}$, we write $x=(x',x_{n})\in \{0,1\}^{n}$. Let $
A_0=\{x'\in\{0,1\}^{n-1}:(x',0)\in A\}$ and $A_1=\{x'\in\{0,1\}^{n-1}:(x',1)\in A\}$. Let $a_0:=\mu^{n-1}(A_0), a_1:=\mu^{n-1}(A_1)$, and note that $a:=\mu^{n}(A)=\frac{a_0+a_1}{2}$. 

\begin{fact}
If $A$ is an increasing set, then $A_0\subseteq A_1$. 
\end{fact}

Define $f_z:A_z\to\mathbb{R}$ by $f_z(x_1,\dots,x_{n-1}):=f(x_1,\dots,x_{n-1},z)$ for $z\in\{0,1\}$.

By expanding the Dirichlet form with respect to the conditional expectation on the $n$-th variable, we can readily express the Dirichlet form of the original function in terms of those of its restricted versions.
\begin{lemma}\label{lem:Dirichlet forms2}
Let $A\subseteq\{0,1\}^n$ be increasing. For any $f:A\to\mathbb{R}$, we have
\begin{equation}
\mathcal{E}_A(f)= \frac{a_1}{2a}\cdot \mathcal{E}_{A_1}(f_1)+\frac{a_0}{2a}\cdot \mathcal{E}_{A_0}(f_0)+\frac{a_0}{4a}\cdot\underset{x'\sim A_0}{\mathbb{E}}\left[(f_0(x')-f_1(x'))^2\right].
\end{equation}
\end{lemma}

\begin{proof}
Since $A$ is increasing, we have $A_0\subseteq A_1$. First, separate the contribution of the last coordinate from that of the first $n-1$ coordinates; this gives
\begin{equation*}
\mathcal{E}_A(f)
=\frac{1}{4}\sum_{i=1}^{n-1}\underset{x\sim A}{\mathbb{E}}\left[(f(x)-f(x^{\oplus i}))^2\cdot\mathbbm{1}_{\{x^{\oplus i}\in A\}}\right]+\frac{1}{4}\underset{x\sim A}{\mathbb{E}}\left[(f(x)-f(x^{\oplus n}))^2\cdot\mathbbm{1}_{\{(x',1-x_n)\in A\}}\right].
\end{equation*}
Conditioning on the value of $x_n$, we further rewrite this as
\begin{equation*}
\begin{split}
\mathcal{E}_A(f)
&=\frac{1}{4}\cdot\frac{|A_1|}{|A|}\sum_{i=1}^{n-1}\underset{x'\sim A_1}{\mathbb{E}}\left[(f_1(x')-f_1(x'^{\oplus i}))^2\cdot\mathbbm{1}_{\{(x'^{\oplus i},1)\in A\}}\right]\\
&\quad+\frac{1}{4}\cdot\frac{|A_0|}{|A|}\sum_{i=1}^{n-1}\underset{x'\sim A_0}{\mathbb{E}}\left[(f_0(x')-f_0(x'^{\oplus i}))^2\cdot\mathbbm{1}_{\{(x'^{\oplus i},0)\in A\}}\right]\\
&\quad+\frac{1}{4}\cdot\frac{|A_1|}{|A|}\underset{x'\sim A_1}{\mathbb{E}}\left[(f_1(x')-f_0(x'))^2\cdot\mathbbm{1}_{\{(x',0)\in A\}}\right]\\
&\quad+\frac{1}{4}\cdot\frac{|A_0|}{|A|}\underset{x'\sim A_0}{\mathbb{E}}\left[(f_0(x')-f_1(x'))^2\cdot\mathbbm{1}_{\{(x',1)\in A\}}\right].
\end{split}
\end{equation*}
Since $A_0\subseteq A_1$,  the third summand above can be rewritten as 
\begin{equation*}
\begin{split}
\frac{|A_1|}{|A|}\underset{x'\sim A_1}{\mathbb{E}}\left[(f_1(x')-f_0(x'))^2\cdot\mathbbm{1}_{\{(x',0)\in A\}}\right]
&=\frac{1}{|A|}\sum_{x'\in A_1}(f_1(x')-f_0(x'))^2\cdot\mathbbm{1}_{\{x'\in A_0\}}\\
&=\frac{|A_0|}{|A|}\underset{x'\sim A_0}{\mathbb{E}}\left[(f_1(x')-f_0(x'))^2\right].
\end{split}
\end{equation*}
Therefore,
\begin{equation*}
\begin{split}
\mathcal{E}_A(f)
&=\frac{a_1}{2a}\cdot\frac{1}{4}\sum_{i=1}^{n-1}\underset{x'\sim A_1}{\mathbb{E}}\left[(f_1(x')-f_1(x'^{\oplus i}))^2\cdot\mathbbm{1}_{\{x'^{\oplus i}\in A_1\}}\right]\\
&\quad+\frac{a_0}{2a}\cdot\frac{1}{4}\sum_{i=1}^{n-1}\underset{x'\sim A_0}{\mathbb{E}}\left[(f_0(x')-f_0(x'^{\oplus i}))^2\cdot\mathbbm{1}_{\{x'^{\oplus i}\in A_0\}}\right]\\
&\quad+\frac{|A_0|}{4|A|}\underset{x'\sim A_0}{\mathbb{E}}\left[(f_1(x')-f_0(x'))^2\right]
+\frac{|A_0|}{4|A|}\underset{x'\sim A_0}{\mathbb{E}}\left[(f_0(x')-f_1(x'))^2\right]\\
&=\frac{a_1}{2a}\cdot \mathcal{E}_{A_1}(f_1)+\frac{a_0}{2a}\cdot \mathcal{E}_{A_0}(f_0)+\frac{a_0}{4a}\cdot\underset{x'\sim A_0}{\mathbb{E}}\left[(f_0(x')-f_1(x'))^2\right].
\end{split}
\end{equation*}
This proves the claim.
\end{proof}

Similarly, by conditioning on the $n$-th coordinate, the variance of the original function can be expressed in terms of the variances of its restricted functions.

\begin{lemma}\label{lem:variance}
For any $f:A\to\mathbb{R}$, we have
\begin{equation}
{\rm Var}_A[f]
=\frac{a_1}{2a}\cdot{\rm Var}_{A_1}[f_1]+\frac{a_0}{2a}\cdot{\rm Var}_{A_0}[f_0]+\frac{a_1a_0}{(a_0+a_1)^2}\cdot\left(\underset{x'\sim A_1}{\mathbb{E}}[f_1(x')]-\underset{x'\sim A_0}{\mathbb{E}}[f_0(x')]\right)^2.
\end{equation}
\end{lemma}

\begin{proof}
Note that
\begin{equation*}
{\rm Var}_A[f]
=\underset{x\sim A}{\mathbb{E}}\left[\,\bigl|f-\underset{x\sim A}{\mathbb{E}}f\bigr|^2\,\right]=\frac{a_1}{2a}\cdot\underset{x'\sim A_1}{\mathbb{E}}\left[\,\bigl|f_1-\underset{x\sim A}{\mathbb{E}}f\bigr|^2\,\right]
+\frac{a_0}{2a}\cdot\underset{x'\sim A_0}{\mathbb{E}}\left[\,\bigl|f_0-\underset{x\sim A}{\mathbb{E}}f\bigr|^2\,\right].
\end{equation*}

Let $\mu_1:=\underset{x'\sim A_1}{\mathbb{E}}[f_1(x')]$ and $\mu_0:=\underset{x'\sim A_0}{\mathbb{E}}[f_0(x')]$. Then
\begin{equation}
\underset{x\sim A}{\mathbb{E}}[f(x)]
=\frac{a_1}{2a}\mu_1+\frac{a_0}{2a}\mu_0.
\end{equation}

Thus, by the identity $\mathbb{E}[(X-b)^2]={\rm Var}[X]+(\mathbb{E}[X]-b)^2$, we have
\begin{equation*}
\underset{x'\sim A_z}{\mathbb{E}}\left[\,\bigl|f_z-\underset{x\sim A}{\mathbb{E}}f\bigr|^2\,\right]
={\rm Var}_{A_z}[f_z]+\frac{a_{1-z}^2}{4a^2}(\mu_1-\mu_0)^2, \quad \forall z\in\{0,1\}.   
\end{equation*}
Therefore,
\begin{equation*}
{\rm Var}_A[f]
=\frac{a_1}{2a}\cdot{\rm Var}_{A_1}[f_1]+\frac{a_0}{2a}\cdot{\rm Var}_{A_0}[f_0]+\frac{a_1a_0}{(a_0+a_1)^2}\cdot(\mu_1-\mu_0)^2.\qedhere
\end{equation*}
\end{proof}

\subsection{Induction part}

Our goal is to prove that for every function $f:A\to \mathbb{R}$,
\begin{equation}\label{ineq:main ineq-c}
   \mathcal{E}_A(f) \ge c\cdot(1-\sqrt{1-a})\cdot{\rm Var}_A[f],
\end{equation}
where $c>0$ is a universal constant (we may take $c=0.5$).

It is easy to verify that the theorem holds for $n=1$. Let $n\ge 2$, and suppose the statement of the theorem holds when $n$ is replaced by $n-1$. Applying Lemmas~\ref{lem:Dirichlet forms2} and~\ref{lem:variance} to split the Dirichlet form and variance of the original function into two restricted functions:
\begin{equation}\label{ineq:Dirichlet and variance}
\begin{split}
\mathcal{E}_A(f) &= \frac{a_1}{2a}\cdot \mathcal{E}_{A_1}(f_1)+\frac{a_0}{2a}\cdot \mathcal{E}_{A_0}(f_0)+\frac{a_0}{4a}\cdot\underset{x'\sim A_0}{\mathbb{E}}\left[(f_0(x')-f_1(x'))^2\right];\\
{\rm Var}_A[f] &=\frac{a_1}{2a}{\rm Var}_{A_1}[f_1]+\frac{a_0}{2a}{\rm Var}_{A_0}[f_0] +\frac{a_1a_0}{(a_0+a_1)^2}(\underset{x'\sim A_1}{\mathbb{E}}[f_1(x')]-\underset{x'\sim A_0}{\mathbb{E}}[f_0(x')])^2.
\end{split}
\end{equation}
By inductive assumption, we have
\begin{equation}\label{ineq:induction hypo}
\begin{split}
&\mathcal{E}_{A_1}(f_1)\ge c\cdot(1-\sqrt{1-a_1}){\rm Var}_{A_1}[f_1],\\
&\mathcal{E}_{A_0}(f_0)\ge c\cdot(1-\sqrt{1-a_0}){\rm Var}_{A_0}[f_0].
\end{split}
\end{equation}
Let $C_1=c\cdot a_1\cdot(\sqrt{1-a}-\sqrt{1-a_1})$ and $C_0=c\cdot a_0\cdot (\sqrt{1-a}-\sqrt{1-a_0})$. Thus, by~\eqref{ineq:Dirichlet and variance} and the inequalities in~\eqref{ineq:induction hypo}, in order to  prove the inequality in~\eqref{ineq:main ineq-c}, it suffices to  show that
\begin{equation}\label{ineq:main five-point inequality}
\begin{split}
    &C_1\cdot{\rm Var}_{A_1}[f_1]+C_0\cdot{\rm Var}_{A_0}[f_0]+\frac{a_0}{2}\cdot\underset{x'\sim A_0}{\mathbb{E}}\left[(f_0(x')-f_1(x'))^2\right]\\
    &\underset{?}{\ge} c\frac{a_1a_0}{a_0+a_1}(1-\sqrt{1-a})(\mu_1-\mu_0)^2.
\end{split}
\end{equation}
\begin{lemma}\label{lem:main reduction}
Assume that $
0<c\le \frac{1}{\sqrt{1-a_0}+\sqrt{1-a_1}}$ (say, $c=\frac{1}{2}$). We may assume $f_0$ is constant on $A_0$, and $f_1$ is constant on both $A_0$ and $A_1\setminus A_0$.  
\end{lemma}
\begin{proof}
Let $\mu_1:=\underset{x'\sim A_1}{\mathbb{E}}[f_1(x')]$, $\alpha:=\mathbb{E}_{x'\sim A_0}[f_1]$, $\beta:=\mathbb{E}_{x'\sim A_1\setminus A_0}[f_1]$ and recall $\mu_0:=\mathbb{E}_{x'\sim A_0}[f_0]$. Note that
\begin{equation} \label{eq:decom A1}
\begin{split}
\underset{x'\sim A_1}{\mathbb{E}}[|f_1(x')-\mu_1|^2]&=\frac{|A_0|}{|A_1|}\cdot\frac{1}{|A_0|}\sum\limits_{x'\in A_0}|f_1(x')-\mu_1|^2\\
&\quad+\frac{|A_1\setminus A_0|}{|A_1|}\cdot\frac{1}{|A_1\setminus A_0|}\sum\limits_{x'\in A_1\setminus A_0}|f_1(x')-\mu_1|^2\\
&=\frac{a_0}{a_1}\underset{x'\sim A_0}{\mathbb{E}}[|f_1-\mu_1|^2]+\frac{a_1-a_0}{a_1}\underset{x'\sim A_1\setminus A_0}{\mathbb{E}}[|f_1-\mu_1|^2],
\end{split}
\end{equation}
and
\begin{equation} \label{eq:decom A0}
 \begin{split}
\underset{x'\sim A_0}{\mathbb{E}}\left[(f_0(x')-f_1(x'))^2\right]&= \underset{x'\sim A_0}{\mathbb{E}}\left[|f_0(x')-\mu_0-(f_1(x')-\mu_1)+\mu_0-\mu_1|^2\right]\\
&=\underset{x'\sim A_0}{\mathbb{E}}\left[|f_0-\mu_0|^2\right]+\underset{x'\sim A_0}{\mathbb{E}}\left[|f_1-\mu_1|^2\right]\\
&\qquad+2\underset{x'\sim A_0}{\mathbb{E}}\left[(f_0-\mu_0)(\mu_1-f_1)\right]\\
&\qquad+(\mu_0-\mu_1)^2+2(\mu_0-\mu_1)(\mu_1-\alpha).
\end{split}
\end{equation}
Combining the decompositions in~\eqref{eq:decom A1} and~\eqref{eq:decom A0}, we obtain that the LHS of~\eqref{ineq:main five-point inequality} is equal to 
\begin{equation}
    \begin{split}
        &\left(C_1\frac{a_0}{a_1}+\frac{a_0}{2}\right)\underset{x'\sim A_0}{\mathbb{E}}[|f_1-\mu_1|^2]+C_1\frac{a_1-a_0}{a_1}\underset{x'\sim A_1\setminus A_0}{\mathbb{E}}[|f_1-\mu_1|^2]\\
        &+\left(C_0+\frac{a_0}{2}\right)\cdot{\rm Var}_{A_0}[f_0]+a_0\cdot \underset{x'\sim A_0}{\mathbb{E}}\left[(f_0-\mu_0)(\mu_1-f_1)\right]\\
        &+\frac{a_0}{2}(\mu_0-\mu_1)^2+a_0(\mu_0-\mu_1)(\mu_1-\alpha).
    \end{split}
\end{equation}
Now we pick up the terms involving expectation or variance on $A_0$: 
\begin{equation}
 \begin{split}
&\left(C_1\frac{a_0}{a_1}+\frac{a_0}{2}\right)\cdot\underset{x'\sim A_0}{\mathbb{E}}[|f_1-\mu_1|^2]
+\left(C_0+\frac{a_0}{2}\right)\cdot{\rm Var}_{A_0}[f_0]+a_0\cdot \underset{x'\sim A_0}{\mathbb{E}}\left[(f_0-\mu_0)(\mu_1-f_1)\right]\\
&=\frac{1}{|A_0|}\sum\limits_{x'\in A_0}\Biggl[
\left(C_1\frac{a_0}{a_1}+\frac{a_0}{2}\right)|f_1(x')-\mu_1|^2
+\left(C_0+\frac{a_0}{2}\right)|f_0(x')-\mu_0|^2\\
&\quad
-a_0(f_0(x')-\mu_0)(f_1(x')-\mu_1)
\Biggr]\\
&=\frac{1}{|A_{0}|}\sum_{x'\in A_{0}}\vec{b}_{x'}^{\top}G\vec{b}_{x'}.
\end{split} 
\end{equation}
where
\[
\vec{b}_{x'}=\begin{bmatrix}f_{0}(x')-\mu_{0}\\
f_{1}(x')-\mu_{1}
\end{bmatrix},\\
G=\begin{bmatrix}a_{0}\cdot\frac{2c(\sqrt{1-a}-\sqrt{1-a_{0}})+1}{2} & -\frac{a_{0}}{2}\\
-\frac{a_{0}}{2} & a_{0}\cdot\frac{2c(\sqrt{1-a}-\sqrt{1-a_{1}})+1}{2}
\end{bmatrix}.\\
\]

We claim that the matrix $G$ is positive semidefinite. Note that $G$ is a symmetric matrix, then we shall use the basic fact that a symmetric matrix $G=\begin{bmatrix}u & v\\
w & q
\end{bmatrix}\in \mathbb{R}^{2\times 2}$ is positive semidefinite if and only if $u\ge 0, q\ge 0$ and $uq-vw\ge 0$. 

Let $s:=\sqrt{1-a_0}$ and $t:=\sqrt{1-a_1}$. Then $0\le t\le s\le 1$ and $1-a=(s^2+t^2)/2$. We first show that the diagonal entries of $G$ are nonnegative. Since $\sqrt{\frac{s^2+t^2}{2}}\ge \frac{s+t}{2}$,
we have
\[
1+2c(\sqrt{1-a}-\sqrt{1-a_0})
\ge 1+2c\left(\frac{s+t}{2}-s\right)
=1-c(s-t)\ge 0,
\]
where the last inequality follows from the assumption $c\le \frac{1}{s+t}$. Also, since $\sqrt{1-a}\ge \sqrt{1-a_1}$, we have
\[
1+2c(\sqrt{1-a}-\sqrt{1-a_1})\ge 1\ge 0.
\]
Thus both diagonal entries of $G$ are nonnegative. It remains to verify that $\det(G)\ge 0$. Note that
\begin{equation}
\begin{split}
&\bigl(1+2c(\sqrt{1-a}-\sqrt{1-a_0})\bigr)\bigl(1+2c(\sqrt{1-a}-\sqrt{1-a_1})\bigr)-1\\
&=\left(1+2c\left(\sqrt{\frac{s^2+t^2}{2}}-s\right)\right)\left(1+2c\left(\sqrt{\frac{s^2+t^2}{2}}-t\right)\right)-1\\
&=2c\left(\sqrt{2}\sqrt{s^2+t^2}-s-t\right)\bigl(1-c(s+t)\bigr)\ge 0,
\end{split}
\end{equation}
where the nonnegativity at the last line follows since $2(s^2+t^2)\ge (s+t)^2$ and $1-c(s+t)\ge 0$. Therefore $\det(G)\ge 0$.
  
Define for $u,v\in\mathbb{R}$ the quadratic function
\[
g(u,v):= C_0\,u^2 + C_1\frac{a_0}{a_1}\,v^2 + \frac{a_0}{2}\,\bigl(u+(\mu_0-\mu_1)-v\bigr)^2.
\]
By the positive semidefiniteness of $G$, we know that the quadratic function $g$ is convex since the Hessian of  $g$ is $G$. Thus by Jensen's inequality, we have
\begin{equation*}
\begin{split}
 &\underset{x'\sim A_0}{\mathbb{E}}[g(f_{0}(x')-\mu_{0},f_{1}(x')-\mu_{1})]\\
 &\ge g \left( \underset{x'\sim A_0}{\mathbb{E}}[f_{0}(x')-\mu_{0}], \underset{x'\sim A_0}{\mathbb{E}}[f_{1}(x')-\mu_{1}]\right)=g(0,\alpha-\mu_1),
\end{split}
\end{equation*}
and
$$
\underset{x'\sim A_1\setminus A_0}{\mathbb{E}}[|f_1(x')-\mu_1|^2]\ge\left(\underset{x'\sim A_1\setminus A_0}{\mathbb{E}}[f_1(x')-\mu_1]\right)^2=(\beta-\mu_1)^2.
$$
Hence replacing $f_{j}$ with their averages on the corresponding subsets does not increase the LHS, and meanwhile, does not change the RHS of~\eqref{ineq:main five-point inequality}.
\end{proof}
\subsection{Five-point Inequality Part}
Invoking Lemma~\ref{lem:main reduction}, we now let $\alpha$ be the value of $f_1$ on $A_0$ and $\beta$ be the value of $f_1$ on $A_1\setminus A_0$. Let $\mu_0$ be the value of $f_0$ on $A_0$. Then
\begin{equation*}
 \begin{split}
    {\rm Var}_{A_1}[f_1]&=\underset{x'\sim A_1}{\mathbb{E}}[f_1^2]-\underset{x'\sim A_1}{\mathbb{E}}[f_1]^2\\
    &=\frac{1}{2^{n-1}\cdot a_1}(a_0\cdot2^{n-1}\cdot \alpha^2+(a_1-a_0)\cdot2^{n-1}\cdot \beta^2)\\
    &\quad-\left(\frac{1}{2^{n-1}\cdot a_1}(a_0\cdot2^{n-1}\cdot \alpha+(a_1-a_0)\cdot2^{n-1}\cdot \beta)\right)^2\\
    &=\frac{a_1a_0-a_0^2}{a_1^2}\cdot(\alpha-\beta)^2.
\end{split} 
\end{equation*}
Since $f_0$ is constant on $A_0$, we have ${\rm Var}_{A_0}[f_0]=0$. Thus to verify~\eqref{ineq:main five-point inequality}, it suffices to prove the following five-point inequality with $c=\frac{1}{2}$. 
\begin{equation}\label{ineq:mian Five-point inequality1}
 \begin{split}
    &c(\sqrt{1-a}-\sqrt{1-a_1})\cdot(a_1-a_0)\cdot(\alpha-\beta)^2+\frac{a_1}{2}\cdot(\mu_0-\alpha)^2\\
   &\ge \frac{c}{a_1+a_0}\cdot(1-\sqrt{1-a})\cdot[a_1\cdot(\beta-\mu_0)+a_0(\alpha-\beta)]^2.  
\end{split} 
\end{equation}
Note that
\begin{equation*}
    \begin{split}
     &\sqrt{1-a}-\sqrt{1-a_1}=\frac{1-a-(1-a_1)}{\sqrt{1-a}+\sqrt{1-a_1}} =\frac{a_1-a_0}{2(\sqrt{1-a}+\sqrt{1-a_1})},\\
     &1-\sqrt{1-a}=\frac{1-(1-a)}{1+\sqrt{1-a}}=\frac{a_1+a_0}{2(1+\sqrt{1-a})},\\
     &a_1\cdot(\beta-\mu_0)+a_0(\alpha-\beta)=(a_0-a_1)\cdot(\alpha-\beta)+ a_1\cdot(\alpha-\mu_0).
    \end{split}
\end{equation*}
After these simple calculations, the inequality in~\eqref{ineq:mian Five-point inequality1} becomes
\begin{equation}
 \begin{split}   
  &c(1-\sqrt{1-a_1})\cdot(a_1-a_0)^2\cdot (\alpha-\beta)^2\\
  &+a_1\cdot(1+\sqrt{1-a}-c\cdot a_1)\cdot(\sqrt{1-a}+\sqrt{1-a_1})\cdot(\alpha-\mu_0)^2\\
  &-2\cdot c\cdot(\sqrt{1-a}+\sqrt{1-a_1})\cdot a_1\cdot(a_0-a_1)\cdot(\alpha-\beta)\cdot(\alpha-\mu_0)\ge 0.
\end{split} 
\end{equation}
If $\alpha=\mu_0$, the inequality is trivially satisfied. Thus, assume  $\alpha\neq\mu_0$ and set $T:=\frac{\alpha-\beta}{\alpha-\mu_0}\cdot (a_1-a_0)$. Then, the LHS above can be rewritten as a quadratic polynomial $AT^2 + BT + C$, whose coefficients depend only on $a_1,a$:
\begin{equation*}
 \begin{split}
&A:=c\cdot (1-\sqrt{1-a_1})\ge 0\\
&B:=2\cdot c\cdot a_1\cdot(\sqrt{1-a}+\sqrt{1-a_1})\\
&C:=a_1\cdot(1+\sqrt{1-a}-c\cdot a_1)\cdot(\sqrt{1-a}+\sqrt{1-a_1}).
\end{split} 
\end{equation*}
To ensure the quadratic is nonnegative for all $T$, it suffices to check that its discriminant satisfies $\Delta=B^2 -4AC\le 0$ with $c=0.5$, which will conclude the proof.
\begin{lemma}
   $B^2-4AC\le 0$. 
\end{lemma}
\begin{proof}
Denote $x=a_0, y=a_1$. Simplifying and rearranging, we need to prove that for any $0\le x\le y\le 1$,
\begin{equation*}
\left(\sqrt{1-\frac{x+y}{2}}+1\right)\cdot(y-2+2\sqrt{1-y})\le 0.
\end{equation*}
It is trivially true as $y-2+2\sqrt{1-y}\le 0,\forall y\in[0,1]$.
\end{proof}

\begin{remark}
The profile $h_c(a) = c(1-\sqrt{1-a})$ with $c=1/2$ in Theorem~\ref{thm:main thm2} and its accompanying proof is a natural choice for ease of analysis, although it is indeed \textit{not} the optimal choice specific to our proof. The optimal choice can be derived as follows. By examining our proof presented above, any choice of $h$ (to replace $h_{1/2}$) allows our proof to hold if and only if it guarantees two conditions: the matrix $G$ is positive semidefinite, and the final five-point inequality holds. The former condition is equivalent to requiring that the diagonal entries of $G$ are nonnegative and that $\det{(G)} \ge 0$ simultaneously. We numerically verify that the choice of $h$ is essentially determined by the constraint $\det{(G)} \ge 0$. Setting the equality to hold and solving this functional equation analytically yields the optimal choice $h^*(a)=-\frac{1}{2}\ln(1-c'a)$, where $c'$ is a constant to be determined. To ensure the final five-point inequality holds, we numerically verify that the optimal choice of $c'$ is approximately 0.8. With this choice of $c'$, $h^*$ is in fact sandwiched between our choice $h_{1/2}$ and Fei and Ferreira Pinto Jr.'s choice $h_{1}$. In other words, our proof, when equipped with any other feasible choice of $h$, cannot yield a bound significantly better than the one presented in Theorem~\ref{thm:main thm2}; furthermore, it cannot yield a bound strictly better than that of Fei and Ferreira Pinto Jr. in~\cite{fei2025spectral}. Thus, our choice is natural in the sense that it  not only simplifies a great deal of complicated calculations in our proof, but also  is ``close'' to optimal. In contrast, a simpler choice $h(a) = \frac{a}{2}$ cannot satisfy the two requirements and thus invalidates the proof, which is reasonable since $h(a) \approx h_{1}(a) > h^*(a)$ when $a$ is sufficiently small.
\end{remark}

\begin{remark}
When restricted to the form $h_c(a)=c\bigl(1-\sqrt{1-a}\,\bigr),$ the optimal choice of $c$ is exactly $c=\frac12$, as we now explain. In Lemma~\ref{lem:main reduction}, the application of Jensen’s inequality requires the matrix $G$ to be positive semidefinite. The proof establishes that this holds if and only if
\[
c\le \frac{1}{\sqrt{1-a_0}+\sqrt{1-a_1}}.
\]
Consequently, the reduction step implies that the optimal choice for $c$ is $\frac12$. Indeed, the final five-point inequality also requires $c\le\frac12$ for validity.	
\end{remark}

\section{Proof of Theorem~\ref{thm:main thm1}}\label{sec:biased}
Our argument follows the induction-by-restrictions framework of Samorodnitsky~\cite{Samorodnitsky17}: we condition on the last coordinate, apply the induction hypothesis to the two restrictions $g_0,g_1$, and then reduce the remaining step to a finite-dimensional inequality in the parameters $a_0,a_1$ (the slice measures) and the averages of $g_0,g_1$.
The main difficulty occurring in the $p$-biased setting is that the resulting ``two-point'' inequality is required to be verified uniformly for all $p\in(0,1)$, which leads to the more involved analytic estimates in subsection~\ref{subsec:two-point-ineq}.

\subsection{Preliminaries}
In this subsection, we outline some notation and terminology, and present some basic lemmas for inductive proof of Theorem~\ref{thm:main thm1}.

Recall that $\mu_p=\mu_p^n$ denote the product $p$-biased measure on $\{0,1\}^n$ defined in~\eqref{def:pbiased}. We consider the space of real valued functions on $\{0,1\}^n$, equipped with the inner product
\begin{equation*}
\tup{f,g}_p=\underset{x\sim\mu_p}{\mathbb{E}}[f(x)g(x)]=\sum\limits_{x\in \{0,1\}^n}f(x)\cdot g(x)\mu_p(x).
\end{equation*}
\begin{definition}
For $q\ge 1$, the $q$-norm w.r.t. the $p$-biased measure of $f:\{0,1\}^n\to\mathbb{R}$ is
\[
\|f\|_{q,\mu_p}:=\left(\underset{x\sim\mu_p}{\mathbb{E}}[|f(x)|^q]\right)^{1/q}.
\]
\end{definition}

Consistent with notations in subsection~\ref{sec:Prethm2}, we write $x=(x',x_{n})\in \{0,1\}^{n}$. Let 
$A_0=\{x'\in\{0,1\}^{n-1}:(x',0)\in A\}$ and $A_1=\{x'\in\{0,1\}^{n-1}:(x',1)\in A\}$.

Let $\mu_p^{n-1}(A_0)=a_0, \mu_p^{n-1}(A_1)=a_1$, and $\mu_p^{n}(A)=a=pa_1+(1-p)a_0$. Define $g_j:(\{0,1\}^{n-1},\mu_p^{n-1})\to\mathbb{R}$ by $g_j(x_1,\dots,x_{n-1}):=g(x_1,\dots,x_{n-1},j)$ for $j\in\{0,1\}$ so that $g_j$ is supported on $A_j$. Throughout the proof we assume that $A$ is increasing (as in the statement of Theorem~\ref{thm:main thm1}), which implies $A_0\subseteq A_1$. 

By applying conditional expectation expansion to $n$-th variable of Dirichlet forms, we can readily reduce the Dirichlet form of the original function to that of the restricted functions.  
\begin{lemma}\label{lem:Dirichlet forms1}
For any $g:\{0,1\}^n\to\mathbb{R}$, we have
\begin{equation}
\mathcal{E}_p(g)= p\mathcal{E}_p(g_1)+(1-p)\mathcal{E}_p(g_0)+\frac{1}{4}\|g_1-g_0\|_{2,\mu_p}^2.
\end{equation}
\end{lemma}
\begin{proof}
We expand the definition of $\mathcal E_p$ and separate the contribution of the last coordinate from that of the first $n-1$ coordinates.
For $x=(x',x_n)\in\{0,1\}^{n-1}\times\{0,1\}$ we have
\[
\sum_{y\sim x}\bigl(g(x)-g(y)\bigr)^2
=\sum_{i=1}^{n-1}\bigl(g(x',x_n)-g(x'^{\oplus i},x_n)\bigr)^2
+\bigl(g(x',x_n)-g(x',1-x_n)\bigr)^2.
\]
Taking expectation with respect to $x'\sim\mu_p^{n-1}$ and $x_n\sim\mu_p$, and recalling the normalization factor $\frac14$ in the definition of $\mathcal E_p$, we obtain
\begin{align*}
\mathcal E_p(g)
&=\frac14\underset{x',x_n}{\mathbb E}\left[\sum_{i=1}^{n-1}\bigl(g(x',x_n)-g(x'^{\oplus i},x_n)\bigr)^2\right]
+\frac14\underset{x',x_n}{\mathbb E}\left[\bigl(g(x',x_n)-g(x',1-x_n)\bigr)^2\right]\\
&=p\,\mathcal E_p(g_1)+(1-p)\,\mathcal E_p(g_0)+\frac14\|g_1-g_0\|_{2,\mu_p}^2,
\end{align*}
as claimed.
\end{proof}

We also need to use some basic properties of the function $f(t):=\frac{\log_p t}{t}$.

\begin{lemma}\label{lem:basic property of basic function}
The function $f(t)$ is decreasing and convex for $t\in[0,1],p\in(0,1)$. Moreover, it satisfies the identity
$$
f(\beta\cdot t)=\frac{\log_p (\beta\cdot t)}{\beta\cdot t}=\frac{1}{\beta}\cdot\frac{\log_p t}{t}+\frac{\log_p \beta}{\beta}\cdot\frac{1}{t}.
$$
\end{lemma}
\begin{proof}
Note that $f'(t)=\frac{1-\log t}{\log p \cdot t^2}\le 0$ and
$f''(t)=\frac{-3+2\log t}{\log p\cdot t^3}\ge 0$ when $t\in[0,1],p\in(0,1)$.
\end{proof}

\subsection{Induction Part}
We may assume without loss of generality that $g\ge 0$, since replacing $g$ by its absolute value preserves the RHS of~\eqref{ineq: Samorodnitsky-type p-biased isoper}, while it only decreases the LHS:
    \begin{equation*}
    \begin{split}
\mathcal{E}_p(g)&=\frac{1}{4}\underset{x\sim\mu_p}{\mathbb{E}}\left[\sum\limits_{y\in\{0,1\}^n:y\sim x}|g(x)-g(y)|^2\right]\\
&\ge \frac{1}{4}\underset{x\sim\mu_p}{\mathbb{E}}\left[\sum\limits_{y\in\{0,1\}^n:y\sim x}(|g(x)|-|g(y)|)^2\right]=\mathcal{E}_p(|g|).
    \end{split}
    \end{equation*}
It is easy to check that the theorem holds for $n=1$. Let $n\ge 2$, and suppose the statement of the theorem holds when $n$ is replaced by $n-1$. Using Lemma~\ref{lem:Dirichlet forms1} to split the Dirichlet form of the original function into those of two restricted functions
\begin{equation}\label{eq:decom Dirichlet in pbiased}
\mathcal{E}_p(g)= p\mathcal{E}_p(g_1)+(1-p)\mathcal{E}_p(g_0)+\frac{1}{4}\|g_1-g_0\|_{2,\mu_p}^2.
\end{equation}
where the Dirichlet forms and the $L^2$-norm at the RHS are defined on the $(n- 1)$-dimensional space. By inductive assumption, we have
\begin{equation}\label{ineq:inductionpbaised}
4p\mathcal{E}_p(g_1)\ge \frac{\mathbb{E}_p[g_1]^2}{a_1}\log_pa_1, \quad
4p\mathcal{E}_p(g_0)\ge \frac{\mathbb{E}_p[g_0]^2}{a_0}\log_pa_0.
\end{equation}
Thus, by utilizing~\eqref{eq:decom Dirichlet in pbiased} and~\eqref{ineq:inductionpbaised}, to prove the desired inequality in Theorem~\ref{thm:main thm1}, it suffices to show that 
\begin{equation}\label{ineq:our main goal ineq}
\begin{split} 
&p\cdot\frac{\mathbb{E}_p[g_1]^2}{a_1}\log_pa_1+(1-p)\cdot\frac{\mathbb{E}_p[g_0]^2}{a_0}\log_pa_0+p\cdot\|g_1-g_0\|_{2,\mu_p}^2\\
&\underset{?}{\ge}\frac{\log_p(pa_1+(1-p)a_0)}{pa_1+(1-p)a_0}\cdot\mathbb{E}_p[g]^2.
\end{split}
\end{equation}

Finally, we may reduce to the case that $g_0$ is constant on $A_0$ and $g_1$ is constant on each of $A_0$ and $A_1\setminus A_0$ (and $g_j\equiv 0$ on $A_j^c$). Define
\[
h_0(x):=\left[\frac{1}{\mu_p^{n-1}(A_0)}\sum_{y\in A_0}g_0(y)\mu_p^{n-1}(y)\right]\mathbbm{1}_{A_0}(x),
\]
and
\begin{equation}
    \begin{split}
h_1(x)&:=\left[\frac{1}{\mu_p^{n-1}(A_0)}\sum_{y\in A_0}g_1(y)\mu_p^{n-1}(y)\right]\mathbbm{1}_{A_0}(x)\\
&+\left[\frac{1}{\mu_p^{n-1}(A_1\setminus A_0)}\sum_{y\in A_1\setminus A_0}g_1(y)\mu_p^{n-1}(y)\right]\mathbbm{1}_{A_1\setminus A_0}(x),   
    \end{split}
\end{equation}
which are respectively constant on $A_0$ and  constant on each of $A_0$ and $A_1\setminus A_0$, as desired. Then $\mathbb E_p[h_0]=\mathbb E_p[g_0]$ and $\mathbb E_p[h_1]=\mathbb E_p[g_1]$, so the RHS in the inequality of \eqref{ineq:our main goal ineq} is unchanged. Moreover,
\begin{equation}
    \begin{split}
\|h_1-h_0\|_{2,\mu_p}^2
&=\frac{1}{\mu_p^{n-1}(A_0)}\left|\sum_{y\in A_0}(g_1(y)-g_0(y))\mu_p^{n-1}(y)\right|^2\\
&\quad+\frac{1}{\mu_p^{n-1}(A_1\setminus A_0)}\left|\sum_{y\in A_1\setminus A_0}g_1(y)\mu_p^{n-1}(y)\right|^2\\
&\le \sum_{y\in A_0}|g_1(y)-g_0(y)|^2\mu_p^{n-1}(y)
+\sum_{y\in A_1\setminus A_0}|g_1(y)|^2\mu_p^{n-1}(y)\\
&=\|g_1-g_0\|_{2,\mu_p}^2,
    \end{split}
\end{equation}
where we have used the Cauchy--Schwarz inequality. Therefore replacing $(g_0,g_1)$ by $(h_0,h_1)$ can only decrease the LHS of the second inequality in \eqref{ineq:our main goal ineq}, while leaving the RHS unchanged. Thus we may assume that $g_0$ is constant on $A_0$ and $g_1$ is constant on both $A_0$ and $A_1\setminus A_0$.

Now let $\alpha$ be the value of $g_1$ on $A_0$, let $\beta$ be the value of $g_1$ on $A_1\setminus A_0$, and let $\mu_0$ be the value of $g_0$ on $A_0$. Then
\begin{equation*}
\begin{split}
&s_1:=\mathbb{E}_p[g_1]=\sum_{x\in A_0}g_1(x)\mu_p^{n-1}(x)
+\sum_{x\in A_1\setminus A_0}g_1(x)\mu_p^{n-1}(x)=\alpha a_0+\beta(a_1-a_0),
\\[1ex]
&s_0:=\mathbb{E}_p[g_0]
=\sum_{x\in A_0}g_0(x)\mu_p^{n-1}(x)
=\mu_0 a_0,
\\[1ex]
&\|g_1-g_0\|_{2,\mu_p}^2
=\sum_{x'\in A_0}(g_1(x')-g_0(x'))^2\mu_p^{n-1}(x')
+\sum_{x'\in A_1\setminus A_0}(g_1(x')-g_0(x'))^2\mu_p^{n-1}(x')\\
&\qquad\quad\quad\quad \ \  =(\alpha-\mu_0)^2 a_0+\beta^2(a_1-a_0).
\end{split} 
\end{equation*}
Then
\begin{equation*}
\mathbb{E}_p[g]=  p\mathbb{E}_{x'\sim \mu_p^{n-1}}[g_1(x')]
+(1-p)\mathbb{E}_{x'\sim \mu_p^{n-1}}[g_0(x')]=ps_1+(1-p)s_0. 
\end{equation*}
With this notation, the inequality to be verified in~\eqref{ineq:our main goal ineq} is
\begin{equation}
\begin{split}
p\cdot\frac{\log_p a_1}{a_1}\cdot s_1^2
+(1-p)\cdot\frac{\log_p a_0}{a_0}\cdot s_0^2
+p\cdot a_0(\alpha-\mu_0)^2
+p\cdot (a_1-a_0)\beta^2\\
\ge \frac{\log_p(pa_1+(1-p)a_0)}{pa_1+(1-p)a_0}\cdot(ps_1+(1-p)s_0)^2.
\end{split}
\end{equation}
Recall that $f(t)=\frac{\log_p t}{t}$ for $t\in(0,1]$. Then our goal is to prove that
\begin{equation}\label{ineq:goal-thm1-2}
\begin{split}
ps_1^2f(a_1)+(1-p)s_0^2f(a_0)
+p a_0(\alpha-\mu_0)^2+p(a_1-a_0)\beta^2\\
\ge (ps_1+(1-p)s_0)^2f(pa_1+(1-p)a_0).
\end{split}
\end{equation}
Expressing $\mu_0$ and $\beta$ as functions of $s_i,a_i$ and of $\alpha$, LHS of~\eqref{ineq:goal-thm1-2} is a quadratic in $\alpha$ with coefficients depending on $s_i$ and $a_i$: LHS of~\eqref{ineq:goal-thm1-2} $=A\alpha^2+B\alpha +C$, where
\begin{equation*}
\begin{split}
&A:=p\left(a_0+\frac{a_0^2}{a_1-a_0}\right)\\
&B:=-2p\left(s_0+\frac{a_0s_1}{a_1-a_0}\right)\\
&C:=ps_1^2f(a_1)+(1-p)s_0^2f(a_0)+p\frac{s_0^2}{a_0}+p\frac{s_1^2}{a_1-a_0}.
\end{split}
\end{equation*}
Minimizing this quadratic polynomial over all possible $\alpha \ge 0$ (with the minimum attained at $\alpha=-\frac{B}{2A}=\frac{s_0a_1-s_0a_0+a_0s_1}{a_1a_0}$), we arrive, after some simple calculations, to the following inequality we need to verify:
\[
ps_1^2f(a_1)+(1-p)s_0^2f(a_0)+p\frac{(s_0-s_1)^2}{a_1}\ge(ps_1+(1-p)s_0)^2f(pa_1+(1-p)a_0).
\]
Next, let $T=\frac{s_0}{s_1}$. The above inequality transforms to a quadratic inequality in $T$:
\begin{equation}
pf(a_1)+(1-p)f(a_0)T^2+p\frac{(1-T)^2}{a_1}\ge(p+(1-p)T)^2f(pa_1+(1-p)a_0).   
\end{equation}
We need to check $h(T):=A'T^2+B'T+C'\ge 0$ with the coefficients
\begin{equation*}
\begin{split}
&A':=(1-p)f(a_0)+\frac{p}{a_1}-(1-p)^2f(pa_1+(1-p)a_0)\\
&B':=-\frac{2p}{a_1}-2p(1-p)f(pa_1+(1-p)a_0)\\
&C':=pf(a_1)+\frac{p}{a_1}-p^2f(pa_1+(1-p)a_0).
\end{split}
\end{equation*}
We shall verify $A'\ge 0$ and $B'^2 -4A'C'\le 0$, which will conclude the proof.

By Lemma~\ref{lem:basic property of basic function}, $f(pa_1+(1-p)a_0))\le f(a_0)$ and
$$
A'\ge \frac{(1-p)a_1f(a_0)+p-(1-p)a_1f(a_0)}{a_1}=\frac{p}{a_1}\ge 0.
$$

Hence it remains to verify the inequality $4A'C'\ge B'^2$, which reduces to a two-point inequality:
\begin{equation}\label{ineq:main two-point inequality}
\begin{split}
&pf(a_1)+(1-p)f(a_0)+(1-p)a_1f(a_0)f(a_1)\\
&\ge (p(1-p)a_1f(a_0)+(1-p)^2a_1f(a_1)+1)f(pa_1+(1-p)a_0).
\end{split}
\end{equation}
\subsection{Two-point Inequality Part}\label{subsec:two-point-ineq}
Renaming the variables $x=a_0$ and $y=a_1$, and recalling the constraints on $a_0$ and $a_1$, it suffices to prove~\eqref{ineq:main two-point inequality}, i.e., \eqref{ineq:main two-point ineq-thm1-2} for $0<x<y\le 1$, since the case $x=y$ is immediate and the case $x=0$ follows by continuity.
\begin{equation}\label{ineq:main two-point ineq-thm1-2}
\begin{split}
&p\cdot f(y)+(1-p)\cdot f(x)+(1-p)\cdot yf(x)f(y)\\
&\ge (p(1-p)yf(x)+(1-p)^2yf(y)+1)\cdot f(py+(1-p)x).
\end{split}
\end{equation}

Rearranging, this is easily seen to be equivalent to
\begin{equation}\label{ineq:main aim}
\Delta_p(x,y)\ge \frac{p(1-p)^2y(f(x)-f(y))^2}{p(1-p)yf(x)+(1-p)^2yf(y)+1},
\end{equation}
where 
$$
\Delta_p(x,y):=pf(y)+(1-p)f(x)-f(py+(1-p)x).
$$
Note that $\Delta_p(x,y)\ge 0$ since $f(t):=\frac{\log_p t}{t}$ is convex.

We now substitute $x=ty$ in $\Delta_p(x,y)$, with $0<t<1$; the endpoint cases $t=0$ and $t=1$ are obtained by continuity (and are immediate), so it is enough to consider the open interval. We then expand $\Delta_p(x,y)$ using Lemma~\ref{lem:basic property of basic function} as follows:
\begin{equation*}
\begin{split}
&\Delta_p(x,y)=\Delta_p(ty,y)=  pf(y)+(1-p)f(t\cdot y)-f((p+(1-p)t)\cdot y)\\
&=\frac{p(1-p)(t-1)^2}{t(p+(1-p)t)}\cdot f(y)+\frac{1}{y}\cdot\left(\frac{(1-p)\log_pt}{t}-\frac{\log_p(p+(1-p)t)}{p+(1-p)t}\right).
\end{split}
\end{equation*}
As to the RHS of \eqref{ineq:main aim}, we have
\begin{equation*}
{\rm RHS}(x,y)={\rm RHS}(ty,y)=\frac{y\left[(\frac{1}{t}-1)f(y)+\frac{1}{y}\cdot\frac{\log_p t}{t}\right]^2}{\left[\frac{1}{t(1-p)}+\frac{1}{p}\right]yf(y)+\frac{\log_p t}{(1-p)t}+\frac{1}{p(1-p)^2}}.
\end{equation*}
Taking $z:=yf(y)$,
\[
\Delta_p(x,y)\ge {\rm RHS} \Leftrightarrow  (Az+B)(Ez+F)\ge (Cz+ D)^2,
\]
where $z=yf(y)=\log_p (y)> 0$ as $0<y<1$, and
\[
\begin{cases}
A=\frac{p(1-p)(t-1)^2}{t(p+(1-p)t)},\\
B=\frac{(1-p)\log_pt}{t}-\frac{\log_p(p+(1-p)t)}{p+(1-p)t},\\
C=\frac{1-t}{t},\\
D=\frac{\log_p t}{t},\\
E=\frac{1}{t(1-p)}+\frac{1}{p}\ge 0,\\
F=\frac{\log_p t}{(1-p)t}+\frac{1}{p(1-p)^2}\ge 0
\end{cases}
\]
depend only on $t$ and $p$.

Observe that
$$
AE=\frac{(t-1)^2}{t^2}=C^2.
$$
Therefore, this reduces to a linear inequality in $z$:
$$
(AF+BE-2CD)\cdot z \ge D^2- BF.
$$
This holds for all non-negative $z$ if and only if $AF+BE\ge 2 CD$ and $BF\ge D^2$.

\begin{lemma}
For any $t\in(0,1)$ and $p\in(0,1)$ we have $AF+BE\ge 2 CD$.
\end{lemma}
\begin{proof}
Simplifying and rearranging terms, the inequality reduces to verifying that
$$
F(t,p):=(1-p)t\log_p t+p(t-1)^2-(p+(1-p)t)\log_p(p+(1-p)t)\ge 0.
$$
We compute the derivatives of $F(t,p)$ with respect to $t$:
\begin{equation*}
\begin{split}
&\frac{\partial}{\partial t} F(t,p)= \frac{1-p}{\ln p}\cdot \ln\left(\frac{t}{p+(1-p)t}\right)+2p\cdot(t-1),\\
&\frac{\partial^2}{\partial t^2} F(t,p)=\frac{p(1-p)}{\ln p}\cdot\frac{1}{t(p+(1-p)t)}+2p,\\
&\frac{\partial^3}{\partial t^3} F(t,p)=\frac{p(1-p)}{\ln(1/p)}\cdot\frac{p+2(1-p)t}{t^2(p+(1-p)t)^2}\ge 0.
\end{split}
\end{equation*}
Hence $\frac{\partial}{\partial t} F(t,p)$ is a convex function of $t$, implying that it can vanish at most twice on the interval $t\in(0,1)$. Note that $\frac{\partial}{\partial t} F(t,p)|_{t\to 0}>0$, $\frac{\partial}{\partial t} F(t,p)|_{t=1}=0$. Therefore, there exists $t_0:=t_0(p)\in(0,1)$ such that $\frac{\partial}{\partial t} F(t,p)|_{t=t_0}=0$ and $\frac{\partial}{\partial t} F(t,p)\ge 0$ for $t\in(0,t_0]$ and $\frac{\partial}{\partial t}F(t,p)$ for $t\in(t_0,1)$. Since $F(0,p)=F(1,p)=0$, the function $F(t,p)$ first increases from $0$ at $t=0$ and then decreases back to $0$ at $t=1$. Consequently, $F(t,p)\ge 0$ for all $t\in(0,1)$, which establishes the desired inequality.
\end{proof}
\begin{lemma}
For any $t\in(0,1)$ and $p\in(0,1)$ we have $BF\ge D^2$.
\end{lemma}
\begin{proof}
Simplifying and rearranging terms, this inequality reduces to showing that
\begin{equation*}
    \begin{split}
        &G(t,p):=(1-p)(p+(1-p)t)\log_pt\\
        &-p(1-p)\log_p(t)\cdot\log_p(p+(1-p)t)-t\cdot\log_p(p+(1-p)t)\ge 0.
    \end{split}
\end{equation*}
Clearly, $G(0,p)=G(1,p)=0$. We aim to show that $G(t,p)$ increases from $0$ at $t=0$, attains its maximum at some interior point, and then decreases back to $0$ at $t=1$; hence $G(t,p)\ge0$ for all $t\in(0,1)$. 

We begin by computing the derivative:
\begin{equation*}
\begin{split}
\frac{\partial}{\partial t} G(t,p)&= \frac{1}{\ln p}\left(\frac{p(1-p)}{t}+(1-p)^2\ln t+(1-p)^2\right)\\
&-\frac{1}{\ln p}\left(\frac{p(1-p)}{\ln p}\cdot\frac{\ln(p+(1-p)t)}{t}\right)\\
&-\frac{1}{\ln p}\left(\frac{p(1-p)^2}{\ln p}\cdot \frac{\ln t}{p+(1-p)t}+\ln(p+(1-p)t)+\frac{(1-p)t}{p+(1-p)t}\right).
\end{split}
\end{equation*}
Define $H(t,p):=t(p+(1-p)t)\cdot\frac{\partial}{\partial t} G(t,p)$. To analyze its sign, we consider successive derivatives of $H$:
\begin{equation*}
\begin{split}
 \frac{\partial}{\partial t} H(t,p)=&\frac{1}{(\ln p)^2}\left[-(p-1)^2p\left(2+\ln t+\ln(p+t-pt)\right)\right]\\
 &+\frac{1}{\ln p}\left[(-2t+p(2t-1))\ln(p+t-pt)\right]\\
&-\frac{1}{\ln p}\left[3(p-1)p(1+p(t-1)-2t)+(p-1)^2(-2t+p(2t-1))\ln t\right];
 \end{split}
\end{equation*}
\begin{equation*}
    \begin{split}
\frac{\partial^2}{\partial t^2} H(t,p)&=\frac{1}{(\ln p)^2}\left(-(p-1)^2p\left(\frac{1}{t}+\frac{1-p}{p+t-pt}\right)\right)\\
&+\frac{1}{\ln p}\left(-3(-2+p)(p-1)p-\frac{(p-1)^2(-2t+p(2t-1))}{t}\right)\\
&+\frac{1}{\ln p}(\frac{(1-p)(-2t+p(2t-1))}{p+t-pt}-2(p-1)^3\ln t+2(p-1)\ln(p+t-pt));        
    \end{split}
\end{equation*}
and
$$
\frac{\partial^3}{\partial t^3} H(t,p)=-\frac{(p-1)^2p}{t^2(p+t-pt)^2(\ln p)^2}\cdot h(t,p),
$$
where $h(t,p):=-2t^2+2pt(-1+2t)+p^2(-1+2t-2t^2)+(-6p(-1+t)t^2+4t^3+p^2(t-1)^2(1+2t))\ln p$. Thus it suffices to consider the sign of function $h(t,p)$ for $t\in(0,1)$ and $p\in (0,1)$. Note that
\begin{equation*}
    \begin{split}
 &\frac{\partial}{\partial t}  h(t,p)=-2(p-1)(-2t+p(-1+2t))+6(p-2)(p(t-1)-t)t\ln p,\\
 &\frac{\partial^2}{\partial t^2}  h(t,p)=-4(p-1)^2+6(p-2)(-2t+p(-1+2t))\ln p,\\
 &\frac{\partial^3}{\partial t^3}  h(t,p)=12(p-2)(p-1)\ln p<0.
    \end{split}
\end{equation*}
Hence $\frac{\partial^2}{\partial t^2}  h(t,p)$ is strictly decreasing on $(0,1)$. Evaluating at $t\to0$, we obtain
$$
\lim_{t\to 0} \frac{\partial^2}{\partial t^2}  h(t,p)=f(p):=-4(p-1)^2-6(p-2)p\ln p.
$$
Note that $f(0)=-4,f(1)=0$, and
\begin{equation*}
    \begin{split}
 &f'(p)=-2(-10+7p+6(p-1)\ln p),\\
 &f''(p)=-26+\frac{12}{p}-12\ln p,\\
 &f'''(p)=-\frac{12(1+p)}{p^2}<0.
    \end{split}
\end{equation*}
Therefore $f''(p)$ has a unique zero $p_0\in(0,1)$ since $\lim_{p\to 0^+}f''(p)=+\infty$ and $f''(1)=-14$. It follows that $f'(p)$ increases on $(0,p_0)$ and decreases on $(p_0,1)$. 

As $\lim_{p\to 0^+}f'(p)=-\infty$ and $f'(1)=6$, there exists a unique point $p_1\in(0,p_0)$ such that $f'(p_1)=0$. Hence $f'(p)<0$ on $(0,p_1)$ and $f'(p)>0$ on $(p_1,1)$, and thus, $f(p)$ first decreases on $(0,p_1)$ and then increases on $(p_1,1)$. Since $f(1)=0$ and $\lim_{p\to 0^+}f(p)=-4$, we conclude that
$$
f(p)=\lim_{t\to 0^+}\frac{\partial^2}{\partial t^2}h(t,p)\le 0
$$
for all $p\in(0,1)$. Consequently, $\frac{\partial^2}{\partial t^2}h(t,p)<0$ for $t\in(0,1)$, and hence $\frac{\partial}{\partial t}h(t,p)$ is strictly decreasing in $t$. Since
$$
\left.\frac{\partial}{\partial t}h(t,p)\right|_{t=0}=2p(p-1)<0,
$$
we have $\frac{\partial}{\partial t}h(t,p)<0$ for all $t\in(0,1)$. Moreover,
$$
h(0,p)=p^2(\ln p-1)<0,
$$
so $h(t,p)<0$ throughout $t\in(0,1)$. Thus, $\frac{\partial^3}{\partial t^3}H(t,p)>0$.

Next, note that $\lim_{t\to0^+}\frac{\partial^2}{\partial t^2} H(t,p)=-\infty$ and $\frac{\partial^2}{\partial t^2} H(t,p)|_{t=1}>0$ with $\frac{\partial^2}{\partial t^2} H(t,p)$ increasing on $t\in(0,1)$. Hence it has a unique zero $t_0\in(0,1)$, and $\frac{\partial}{\partial t} H(t,p)$ decreases on $(0,t_0)$ and increases on $(t_0,1)$. Moreover, $\lim_{t\to0^+}\frac{\partial}{\partial t}H(t,p)>0$ and $\frac{\partial}{\partial t}H(t,p)|_{t=1}>0$. Since $
H(1,p)=\lim_{t\to0^+}H(t,p)=0$, we have
\[
\int_0^1 \frac{\partial}{\partial t}H(t,p)\,dt=0.
\]
Therefore $\frac{\partial}{\partial t}H(t,p)$ must be negative at some point of $(0,1)$. As $\frac{\partial}{\partial t}H(t,p)$ attains its unique minimum at $t_0$, it follows that $\frac{\partial}{\partial t}H(t,p)|_{t=t_0}<0$. Hence there exist $0<t_1<t_2<1$ such that $H(t,p)$ first increases on $(0,t_1)$, then decreases on $[t_1,t_2]$, and finally increases again on $(t_2,1)$. Since $\lim_{t\to0^+}H(t,p)=H(1,p)=0$, there exists $t'\in[0,1]$ such that $H(t,p)>0$ for $t\in(0,t')$ and $H(t,p)\le 0$ for $t\in[t',1)$. It follows that $G(t,p)$ increases on $t\in(0,t')$ and then decreases on $t\in[t',1)$. Since $\lim_{t\to0^+}G(t,p)=G(1,p)=0$, we conclude that $G(t,p)\ge 0$ for all $t\in[0,1]$.
\end{proof}

\subsection{Equality case in Theorem~\ref{thm:main thm1}}

For completeness, we explain how the ``only if'' direction in the equality statement follows from the proof. We argue by induction on $n$. The case $n=1$ is immediate.

Assume now that equality holds in~\eqref{ineq: Samorodnitsky-type p-biased isoper} for some nonzero $g$ supported on an increasing set $A\subseteq\{0,1\}^n$. Since replacing $g$ by $|g|$ leaves the RHS unchanged and can only decrease the LHS, equality for $g$ implies equality for $|g|$ as well. Thus we may assume $g\ge 0$.

Recall the inductive decomposition
\begin{equation*}
\begin{split}
4p\mathcal{E}_p(g)
&\ge p\cdot\frac{\mathbb{E}_p[g_1]^2}{a_1}\log_pa_1
+(1-p)\cdot\frac{\mathbb{E}_p[g_0]^2}{a_0}\log_pa_0
+p\cdot\|g_1-g_0\|_{2,\mu_p}^2\\
&\ge \frac{\log_p(pa_1+(1-p)a_0)}{pa_1+(1-p)a_0}\cdot\mathbb{E}_p[g]^2.
\end{split}
\end{equation*}
Hence equality in~\eqref{ineq: Samorodnitsky-type p-biased isoper} forces equality at every subsequent step of the proof.

First, equality in the averaging reduction implies that $g_0$ is constant on $A_0$, and that $g_1$ is constant on each of $A_0$ and $A_1\setminus A_0$.

Next, equality in the final two-point inequality can occur only at the boundary cases. Indeed, in subsection~\ref{subsec:two-point-ineq} we prove that the relevant one-variable function is strictly positive for $0<t<1$, so equality is possible only when $t=0$ or $t=1$, i.e., only when either $a_0=0$ or $a_0=a_1$.

We now consider these two cases.

\smallskip
\noindent\textbf{Case 1: $a_0=0$.}
Then $A_0=\varnothing$, so
\[
A=\{(x',1):x'\in A_1\}=A_1\times\{1\}.
\]
Moreover, $g_0\equiv 0$, and equality in the first inductive inequality yields equality for $g_1$ on $A_1$. By the induction hypothesis, $A_1$ is an increasing subcube and $g_1=c\,\mathbbm{1}_{A_1}$ for some constant $c\ge 0$. Therefore $A=A_1\times\{1\}$ is an increasing subcube and $g=c\cdot\mathbbm{1}_A$.

\smallskip
\noindent\textbf{Case 2: $a_0=a_1$.}
Since $A_0\subseteq A_1$ and $\mu_p^{n-1}(A_0)=\mu_p^{n-1}(A_1)$, we must have $A_0=A_1$. Hence
\[
A=A_0\times\{0,1\}.
\]
Equality in the first inductive inequality now yields equality for both $g_0$ on $A_0$ and $g_1$ on $A_1=A_0$. By the induction hypothesis,
\[
A_0 \text{ is an increasing subcube},\qquad
g_0=c_0\,\mathbbm{1}_{A_0},\qquad
g_1=c_1\,\mathbbm{1}_{A_0}
\]
for some constants $c_0,c_1\ge 0$.

It remains to show that $c_0=c_1$. Since $a_0=a_1$, the inequality to be verified in the induction step becomes
\begin{equation*}
p\,a_0\log_p(a_0)\,c_1^2+(1-p)\,a_0\log_p(a_0)\,c_0^2+p\,a_0(c_1-c_0)^2
\ge a_0\log_p(a_0)\,(pc_1+(1-p)c_0)^2.
\end{equation*}
Subtracting the right-hand side from the left-hand side, we obtain
\begin{equation*}
a_0\,p\bigl(1+(1-p)\log_p(a_0)\bigr)(c_1-c_0)^2\ge 0.
\end{equation*}
Since equality holds and $1+(1-p)\log_p(a_0)>0$, it follows that $c_1=c_0$. Therefore 
\[
g=c_0\mathbbm{1}_A,
\]
and $A=A_0\times\{0,1\}$ is an increasing subcube.

In either case, equality implies that $A$ is an increasing subcube and $g$ is a scalar multiple of $\mathbbm{1}_A$. The converse is immediate from a direct substitution.

\section{Proof of Theorem~\ref{thm:mean-first-exit-time}}\label{sec:mean first exit time}

We now derive Theorem~\ref{thm:mean-first-exit-time} from the $p$-biased functional inequality~\eqref{ineq: Samorodnitsky-type p-biased isoper}.

\begin{proof}[Proof of Theorem~\ref{thm:mean-first-exit-time}]
Let $\mathbf P=\mathbf P_p^{\mathrm{Gl}}$, and let $P_A$ be the killed kernel on $A$, defined by $P_A(x,y)=\mathbf P(x,y)$ for $x,y\in A$. Since $A\subsetneq\{0,1\}^n$, its complement is non-empty. Moreover, for $0<p<1$, the Glauber dynamics $\mathbf P$ is irreducible on the finite cube. Hence the first exit time $\tau$ from $A$ has finite expectation from every starting state in $A$. Define $u(x):=\E_x[\tau]$ for $x\in A$, and extend $u$ by $0$ on $A^c$. By the Markov property, for each $x\in A$,
$$
u(x)=1+\sum_{y\in A}P_A(x,y)u(y),
$$
that is, $(I-P_A)u(x)=1$ for any $x\in A$. Taking the $\mu_p$-inner product with $u$, and using that $u=0$ on $A^c$, we obtain
\begin{equation}\label{eq:inner-poisson}
\langle u,(I-P_A)u\rangle_p=\E_{\mu_p}[u]=\mu_p(A)\E_\nu[\tau].
\end{equation}

Furthermore, since $u\equiv 0$ on $A^c$, we have $(P_Au)(x)=(\mathbf Pu)(x)$ for $x\in A$, and therefore
\begin{equation}\label{eq:PA-to-P}
\langle u,(I-P_A)u\rangle_p=\langle u,(I-\mathbf P)u\rangle_p.
\end{equation}

Recall that
$$
\mathcal E_p(g):=\frac14\sum_{x\in\{0,1\}^n}\mu_p(x)\sum_{y\sim x}\bigl(g(x)-g(y)\bigr)^2.
$$
Since $\mathbf P$ is reversible with respect to $\mu_p$, we have
$$
\langle g,(I-\mathbf P)g\rangle_p
=\frac12\sum_{x,y}\mu_p(x)\mathbf P(x,y)\bigl(g(x)-g(y)\bigr)^2.
$$
For an edge $\{x,x^{\oplus i}\}$, a direct computation gives
$$
\mu_p(x)\mathbf P(x,x^{\oplus i})
=\mu_p(x^{\oplus i})\mathbf P(x^{\oplus i},x)
=\frac{p(1-p)}{n}\mu_p(x_{\backslash i}).
$$
Summing the two directed contributions across each undirected edge, we obtain
\begin{equation}\label{eq:dirichlet-identity}
\langle g,(I-\mathbf P)g\rangle_p=\frac{4p(1-p)}{n}\cdot\mathcal E_p(g).
\end{equation}
Applying \eqref{eq:dirichlet-identity} to $g=u$, and combining \eqref{eq:inner-poisson} and \eqref{eq:PA-to-P}, yields
\begin{equation}\label{eq:Ep-u-final}
\mathcal E_p(u)=\frac{n}{4p(1-p)}\cdot\mu_p(A)\cdot\E_\nu[\tau].
\end{equation}

Since $u\ge 0$ and is supported on $A$, Theorem~\ref{thm:main thm1} applied to $g=u$ gives
$$
4p\cdot\mathcal E_p(u)\ge \frac{\E_{\mu_p}[u]^2}{\mu_p(A)}\log_p\mu_p(A).
$$
Because $A$ is non-empty and proper, we have $0<\mu_p(A)<1$, and hence $\log_p(\mu_p(A))>0$. Substituting \eqref{eq:Ep-u-final} and $\E_{\mu_p}[u]=\mu_p(A)\E_\nu[\tau]$, we obtain
$$
\frac{n}{1-p}\cdot\mu_p(A)\cdot\E_\nu[\tau]
\ge
\mu_p(A)\cdot\E_\nu[\tau]^2\cdot\log_p\mu_p(A).
$$
Since $\mu_p(A)\E_\nu[\tau]>0$, dividing both sides by this quantity gives
$$
\E_\nu[\tau]\le \frac{n}{(1-p)\log_p(\mu_p(A))}.
$$

If equality holds, then equality must hold in Theorem~\ref{thm:main thm1} for $g=u$. By the equality characterization in Theorem~\ref{thm:main thm1}, it follows that $A$ is an increasing subcube and that $u$ is a scalar multiple of $\mathbbm{1}_A$. Conversely, equality holds for increasing subcubes by Example~\ref{exam:subcube-exit}.
\end{proof}

\section{Conclusion and open problems}\label{sec:discuss}
A natural next target beyond the $O(n^2)$ mixing bound for constant-density increasing sets is the full Ding--Mossel conjecture asking for an $O(n\log n)$ mixing time for the censored walk when $\mu(A)\ge\varepsilon>0$ (\cite[Question 1.1]{DM2014}). As observed in
\cite{fei2025spectral}, a route to such an $O(n\log n)$ bound is to establish an appropriate (standard or modified) log-Sobolev inequality for the censored chain on $A$. Following the general relations between log-Sobolev and mixing time bound (e.g.,\cite{BT2006logsobolev,DS1996logsobolev}), it is natural to ask what the \emph{correct} dependence on the density $\mu(A)$ should be in any valid log–Sobolev (or modified log-Sobolev) inequality for increasing sets.

Our inductive two-point machinery is well suited for quadratic ($L^2$) control: it converts high-dimensional Dirichlet forms and variances into lower-dimensional analogues plus a finite list of analytic inequalities (e.g., two-point inequalities) that one can verify exactly. However, in log-Sobolev inequalities, the relative entropy is a nonlinear, non-quadratic functional—this introduces challenges when applying the induction method to control it. Hence,  to establish  log-Sobolev inequalities in non-product spaces, further developing the induction method or new methods are needed. 

\paragraph{Note added.}
After completing this manuscript, Santos, Tripathi, and Youssef~\cite{STY2026} studied logarithmic Sobolev inequalities for the censored walk under random monotone censoring. Writing $C_{\mathrm{LS}}(A)$ for the best constant in the corresponding standard logarithmic Sobolev inequality, they proved that, in the Markov-chain normalization, $C_{\mathrm{LS}}(A_n)\asymp n$ with high probability when $A_n$ is chosen uniformly among all increasing subsets of $\{0,1\}^n$. Consequently, the Ding--Mossel $O(n\log n)$ mixing bound holds for a typical increasing set. They also showed that for every fixed $\alpha\in(0,1/2]$,
\[
\sup_{\substack{A\subseteq\{0,1\}^n\ \mathrm{increasing}\\|A|\ge \alpha 2^n}}C_{\mathrm{LS}}(A)\asymp_\alpha n^2 .
\]
Equivalently, in the normalization of $E_A$ used in $\eqref{eq:restricted-Dirichlet-form}$, a standard logarithmic Sobolev inequality with a constant of order one cannot hold uniformly over all dense increasing sets. Thus, while the full Ding--Mossel conjecture remains open, its proof cannot proceed solely through such a uniform standard log-Sobolev inequality.
\paragraph{Acknowledgments.}
We thank the anonymous referees for their careful reading and helpful suggestions. In particular, we are very grateful to one referee for suggesting us to consider the interpretation of our inequality in Theorem~\ref{thm:main thm1} based on Glauber dynamics, which yields Theorem~\ref{thm:mean-first-exit-time}. This interpretation is more natural than the original $p$-biased random walk considered in the initial version of this paper. 

\bibliographystyle{plain}
\bibliography{references}

@Misc{durcik2024sharpisoperimetricinequalitieshamming,
title={{Sharp isoperimetric inequalities on the {H}amming cube near the critical exponent}},
  author={Durcik, Polona and Ivanisvili, Paata and Roos, Joris},
eprint = {2407.12674},
year={2024},
month = Jul,
eprintclass = {math.CA}
}

@Misc{guruswami2016rapidly,
  title={Rapidly mixing markov chains: a comparison of techniques (a survey)},
  author={Guruswami, Venkatesan},
  eprint = {1603.01512},
  year={2016},
  month = Mar,
eprintclass = {cs.DS}
}

@article {Harper64,
    AUTHOR = {Harper, L. H.},
     TITLE = {Optimal assignments of numbers to vertices},
   JOURNAL = {J. Soc. Indust. Appl. Math.},
  FJOURNAL = {Journal of the Society for Industrial and Applied Mathematics},
    VOLUME = {12},
      YEAR = {1964},
     PAGES = {131--135},
      ISSN = {0368-4245},
   MRCLASS = {05.24},
  MRNUMBER = {162737},
MRREVIEWER = {W. W. Peterson},
}

@article {Hart76,
    AUTHOR = {Hart, Sergiu},
     TITLE = {A note on the edges of the {$n$}-cube},
   JOURNAL = {Discrete Math.},
  FJOURNAL = {Discrete Mathematics},
    VOLUME = {14},
      YEAR = {1976},
    NUMBER = {2},
     PAGES = {157--163},
      ISSN = {0012-365X},
   MRCLASS = {05A99 (90D10)},
  MRNUMBER = {396293},
MRREVIEWER = {G. L. Alexanderson},
       DOI = {10.1016/0012-365X(76)90058-3},
}

@article {Lindsey64,
    AUTHOR = {Lindsey, II, John H.},
     TITLE = {Assignment of numbers to vertices},
   JOURNAL = {Amer. Math. Monthly},
  FJOURNAL = {American Mathematical Monthly},
    VOLUME = {71},
      YEAR = {1964},
     PAGES = {508--516},
      ISSN = {0002-9890},
   MRCLASS = {05.35},
  MRNUMBER = {168489},
MRREVIEWER = {N. S. Mendelsohn},
       DOI = {10.2307/2312587},
}

@article {Bernstein67,
    AUTHOR = {Bernstein, A. J.},
     TITLE = {Maximally connected arrays on the {$n$}-cube},
   JOURNAL = {SIAM J. Appl. Math.},
  FJOURNAL = {SIAM Journal on Applied Mathematics},
    VOLUME = {15},
      YEAR = {1967},
     PAGES = {1485--1489},
      ISSN = {0036-1399},
   MRCLASS = {05.35},
  MRNUMBER = {223260},
MRREVIEWER = {D. W. Walkup},
       DOI = {10.1137/0115129},
}

@article {EKL2019biased,
    AUTHOR = {Ellis, David and Keller, Nathan and Lifshitz, Noam},
     TITLE = {On a biased edge isoperimetric inequality for the discrete
              cube},
   JOURNAL = {J. Combin. Theory Ser. A},
  FJOURNAL = {Journal of Combinatorial Theory. Series A},
    VOLUME = {163},
      YEAR = {2019},
     PAGES = {118--162},
      ISSN = {0097-3165},
   MRCLASS = {05C75 (06E30)},
  MRNUMBER = {3884706},
MRREVIEWER = {Hua Wang},
       DOI = {10.1016/j.jcta.2018.12.001},
}

@article {Samorodnitsky17,
    AUTHOR = {Samorodnitsky, Alex},
     TITLE = {An inequality for functions on the {H}amming cube},
   JOURNAL = {Combin. Probab. Comput.},
  FJOURNAL = {Combinatorics, Probability and Computing},
    VOLUME = {26},
      YEAR = {2017},
    NUMBER = {3},
     PAGES = {468--480},
      ISSN = {0963-5483},
   MRCLASS = {05D05 (60C05)},
  MRNUMBER = {3628914},
 DOI = {10.1017/S0963548316000432},
}

@article {Kahnkalai07,
    AUTHOR = {Kahn, Jeff and Kalai, Gil},
     TITLE = {Thresholds and expectation thresholds},
   JOURNAL = {Combin. Probab. Comput.},
  FJOURNAL = {Combinatorics, Probability and Computing},
    VOLUME = {16},
      YEAR = {2007},
    NUMBER = {3},
     PAGES = {495--502},
      ISSN = {0963-5483},
   MRCLASS = {60C05 (05C80)},
  MRNUMBER = {2312440},
MRREVIEWER = {Amites Sarkar},
       DOI = {10.1017/S0963548307008474},
}

@article {MT2006,
    AUTHOR = {Montenegro, Ravi and Tetali, Prasad},
     TITLE = {Mathematical aspects of mixing times in {M}arkov chains},
   JOURNAL = {Found. Trends Theor. Comput. Sci.},
  FJOURNAL = {Foundations and Trends${}^\circledR$ in Theoretical Computer
              Science},
    VOLUME = {1},
      YEAR = {2006},
    NUMBER = {3},
     PAGES = {x+121},
      ISSN = {1551-305X},
   MRCLASS = {68Q25 (60J10 60J22 90C40)},
  MRNUMBER = {2341319},
       DOI = {10.1561/0400000003},
}

@book {Ryanbook2014,
    AUTHOR = {O'Donnell, Ryan},
     TITLE = {Analysis of {B}oolean functions},
 PUBLISHER = {Cambridge University Press, New York},
      YEAR = {2014},
     PAGES = {xx+423},
      ISBN = {978-1-107-03832-5},
   MRCLASS = {42-02 (06E30 42A16 68-02 68Q15 68Q30 91B14 94D05)},
  MRNUMBER = {3443800},
MRREVIEWER = {Martin C. Cooper},
       DOI = {10.1017/CBO9781139814782},
}

@Misc{yu2025averagedistancelevel1fourier,
  title={{On Average Distance, {L}evel-1 {F}ourier Weight, and {C}hang's Lemma}},
  author={Yu, Lei},
  eprint = {2504.02593},
  year = 2014,
  month = Apr,
  eprintclass = {math.CO}
}

@Misc{beltran2023sharpisoperimetricinequalitieshypercube,
  title={On sharp isoperimetric inequalities on the hypercube},
  author={Beltran, David and Ivanisvili, Paata and Madrid, Jos{\'e}},
  eprint = {2303.06738},
  year={2023},
  month = Mar,
  eprintclass = {math.CA}
}

@Misc{alpay2025lowerboundsdyadicsquare,
  title={Lower Bounds for Dyadic Square Functions of indicator functions of sets},
  author={Alpay, Natanael and Ivanisvili, Paata},
  eprint = {2502.16045},
  year={2025},
  month = Feb,
  eprintclass = {math.CA}
}

@article {FS2007CPCKKL,
    AUTHOR = {Falik, Dvir and Samorodnitsky, Alex},
     TITLE = {Edge-isoperimetric inequalities and influences},
   JOURNAL = {Combin. Probab. Comput.},
  FJOURNAL = {Combinatorics, Probability and Computing},
    VOLUME = {16},
      YEAR = {2007},
    NUMBER = {5},
     PAGES = {693--712},
      ISSN = {0963-5483},
   MRCLASS = {05C12},
  MRNUMBER = {2346808},
MRREVIEWER = {Jos\'{e} C\'{a}ceres},
       DOI = {10.1017/S0963548306008340},
}

@article {KLMM2024globalhyper2024,
    AUTHOR = {Keevash, Peter and Lifshitz, Noam and Long, Eoin and Minzer,
              Dor},
     TITLE = {Hypercontractivity for global functions and sharp thresholds},
   JOURNAL = {J. Amer. Math. Soc.},
  FJOURNAL = {Journal of the American Mathematical Society},
    VOLUME = {37},
      YEAR = {2024},
    NUMBER = {1},
     PAGES = {245--279},
      ISSN = {0894-0347},
   MRCLASS = {06E30 (94D10)},
  MRNUMBER = {4654613},
MRREVIEWER = {I. Villa},
       DOI = {10.1090/jams/1027},
}

@article {DS2005Annal,
    AUTHOR = {Dinur, Irit and Safra, Samuel},
     TITLE = {On the hardness of approximating minimum vertex cover},
   JOURNAL = {Ann. of Math. (2)},
  FJOURNAL = {Annals of Mathematics. Second Series},
    VOLUME = {162},
      YEAR = {2005},
    NUMBER = {1},
     PAGES = {439--485},
      ISSN = {0003-486X},
   MRCLASS = {68Q17 (68Q15)},
  MRNUMBER = {2178966},
MRREVIEWER = {Johan H\aa stad},
       DOI = {10.4007/annals.2005.162.439},
}

@article {Talagrand1994randomgraph,
    AUTHOR = {Talagrand, Michel},
     TITLE = {On {R}usso's approximate zero-one law},
   JOURNAL = {Ann. Probab.},
  FJOURNAL = {The Annals of Probability},
    VOLUME = {22},
      YEAR = {1994},
    NUMBER = {3},
     PAGES = {1576--1587},
      ISSN = {0091-1798},
   MRCLASS = {28A35 (60K35)},
  MRNUMBER = {1303654},
       URL =
              {http://links.jstor.org/sici?sici=0091-1798(199407)22:3<1576:ORAZL>2.0.CO;2-S&origin=MSN},
}

@InProceedings{fei2025spectral,
  author =	{Fei, Yumou and Ferreira Pinto Jr., Renato},
  title =	{{On the Spectral Expansion of Monotone Subsets of the Hypercube}},
  booktitle =	{Approximation, Randomization, and Combinatorial Optimization. Algorithms and Techniques (APPROX/RANDOM 2025)},
  pages =	{42:1--42:24},
  series =	{Leibniz International Proceedings in Informatics (LIPIcs)},
  ISBN =	{978-3-95977-397-3},
  ISSN =	{1868-8969},
  year =	{2025},
  volume =	{353},
  editor =	{Ene, Alina and Chattopadhyay, Eshan},
  publisher =	{Schloss Dagstuhl -- Leibniz-Zentrum f{\"u}r Informatik},
  address =	{Dagstuhl, Germany},
  doi =		{10.4230/LIPIcs.APPROX/RANDOM.2025.42},
}

@book {LP2017,
    AUTHOR = {Levin, David A. and Peres, Yuval},
     TITLE = {Markov chains and mixing times},
      NOTE = {Second edition of [ MR2466937],
              With contributions by Elizabeth L. Wilmer,
              With a chapter on ``Coupling from the past'' by James G. Propp
              and David B. Wilson},
 PUBLISHER = {American Mathematical Society, Providence, RI},
      YEAR = {2017},
     PAGES = {xvi+447},
      ISBN = {978-1-4704-2962-1},
   MRCLASS = {60J10 (60-01 60B15 60C05 60J27 60K35 68U20 82C22)},
   MRCLASS = {60J10 (60-01 60B15 60C05 60J27 60K35 68U20 82C22)},
  MRNUMBER = {3726904},
       DOI = {10.1090/mbk/107},
}

@article {DM2014,
    AUTHOR = {Ding, Jian and Mossel, Elchanan},
     TITLE = {Mixing under monotone censoring},
   JOURNAL = {Electron. Commun. Probab.},
  FJOURNAL = {Electronic Communications in Probability},
    VOLUME = {19},
      YEAR = {2014},
     PAGES = {no. 46, 6},
   MRCLASS = {60J10 (60K35 82B20)},
  MRNUMBER = {3233208},
       DOI = {10.1214/ECP.v19-3157},
}

@article {BKS1999,
    AUTHOR = {Benjamini, Itai and Kalai, Gil and Schramm, Oded},
     TITLE = {Noise sensitivity of {B}oolean functions and applications to percolation},
   JOURNAL = {Inst. Hautes \'{E}tudes Sci. Publ. Math.},
  FJOURNAL = {Institut des Hautes \'{E}tudes Scientifiques. Publications
              Math\'{e}matiques},
    NUMBER = {90},
      YEAR = {1999},
     PAGES = {5--43 (2001)},
      ISSN = {0073-8301},
   MRCLASS = {60B15 (60K35 68Q15 82B43 94C10)},
  MRNUMBER = {1813223},
MRREVIEWER = {H. Kesten},
       URL = {http://www.numdam.org/item?id=PMIHES_1999__90__5_0},
}

@Misc{BIMP2025young,
  title={{Optimal {Y}oung's convolutions inequality and its reverse form on the hypercube}},
  author={Beltran, David and Ivanisvili, Paata and Madrid, Jos{\'e} and Patil, Lekha},
  eprint = {2507.06115},
  year = 2025,
  month = Jul,
  eprintclass = {math.CA}
}

@Misc{CTKS2025inequalities,
  title={{Inequalities in {F}ourier analysis on binary cubes}},
  author={Crmari{\'c}, Ton{\'c}i and Kova{\v{c}}, Vjekoslav and Shiraki, Shobu},
   eprint = {2507.01359},
  year = 2025,
  month = Jul,
  eprintclass = {math.CA}
}

@article {BT2006logsobolev,
    AUTHOR = {Bobkov, Sergey G. and Tetali, Prasad},
     TITLE = {Modified logarithmic {S}obolev inequalities in discrete
              settings},
   JOURNAL = {J. Theoret. Probab.},
  FJOURNAL = {Journal of Theoretical Probability},
    VOLUME = {19},
      YEAR = {2006},
    NUMBER = {2},
     PAGES = {289--336},
      ISSN = {0894-9840},
   MRCLASS = {60J10 (46E35 47D07 60E15 60J35)},
  MRNUMBER = {2283379},
MRREVIEWER = {Ivan Gentil},
       DOI = {10.1007/s10959-006-0016-3},
}

@article {DS1996logsobolev,
    AUTHOR = {Diaconis, P. and Saloff-Coste, L.},
     TITLE = {Logarithmic {S}obolev inequalities for finite {M}arkov chains},
   JOURNAL = {Ann. Appl. Probab.},
  FJOURNAL = {The Annals of Applied Probability},
    VOLUME = {6},
      YEAR = {1996},
    NUMBER = {3},
     PAGES = {695--750},
      ISSN = {1050-5164},
   MRCLASS = {60J05 (15A60 47D07 60F05 60J10)},
  MRNUMBER = {1410112},
MRREVIEWER = {Laurent Miclo},
       DOI = {10.1214/aoap/1034968224},
}

@article {Eldan2025Isoperimetric,
    AUTHOR = {Eldan, Ronen and Kindler, Guy and Lifshitz, Noam and Minzer,
              Dor},
     TITLE = {Isoperimetric inequalities made simpler},
   JOURNAL = {Discrete Anal.},
  FJOURNAL = {Discrete Analysis},
      YEAR = {2025},
     PAGES = {Paper No. 7, 23},
   MRCLASS = {94D10 (06E30)},
  MRNUMBER = {4946554},
}

@article {Ivanisvili2024KKL,
    AUTHOR = {Ivanisvili, Paata and Stone, Yonathan},
     TITLE = {The {KKL} inequality and {R}ademacher type 2},
   JOURNAL = {Discrete Anal.},
  FJOURNAL = {Discrete Analysis},
      YEAR = {2024},
     PAGES = {Paper No. 2, 14},
   MRCLASS = {46B09 (60E15)},
  MRNUMBER = {4741754},
MRREVIEWER = {Pierre Portal},
}

@article {ivanisvili2020rademacher,
    AUTHOR = {Ivanisvili, Paata and van Handel, Ramon and Volberg,
              Alexander},
     TITLE = {Rademacher type and {E}nflo type coincide},
   JOURNAL = {Ann. of Math. (2)},
  FJOURNAL = {Annals of Mathematics. Second Series},
    VOLUME = {192},
      YEAR = {2020},
    NUMBER = {2},
     PAGES = {665--678},
      ISSN = {0003-486X},
   MRCLASS = {46B09 (46B07 60E15)},
  MRNUMBER = {4151086},
MRREVIEWER = {El\'{o}i M. Galego},
       DOI = {10.4007/annals.2020.192.2.8},
}

@article {gross1975logarithmic,
    AUTHOR = {Gross, Leonard},
     TITLE = {Logarithmic {S}obolev inequalities},
   JOURNAL = {Amer. J. Math.},
  FJOURNAL = {American Journal of Mathematics},
    VOLUME = {97},
      YEAR = {1975},
    NUMBER = {4},
     PAGES = {1061--1083},
      ISSN = {0002-9327},
   MRCLASS = {46E35 (81.47)},
  MRNUMBER = {420249},
MRREVIEWER = {R. H\o egh-Krohn},
       DOI = {10.2307/2373688},
}

@incollection {Jer1998,
    AUTHOR = {Jerrum, Mark},
     TITLE = {Mathematical foundations of the {M}arkov chain {M}onte {C}arlo
              method},
 BOOKTITLE = {Probabilistic methods for algorithmic discrete mathematics},
    SERIES = {Algorithms Combin.},
    VOLUME = {16},
     PAGES = {116--165},
 PUBLISHER = {Springer, Berlin},
      YEAR = {1998},
   MRCLASS = {60J10 (05C80 54C05 68Q25 68U20)},
  MRNUMBER = {1678570},
MRREVIEWER = {James Allen Fill},
       DOI = {10.1007/978-3-662-12788-9\_4},
}

@Misc{STY2026,
  title={{Log-Sobolev under random monotone censoring
}},
  author={{Patrick Oliveira Santos, Raghavendra Tripathi, Pierre Youssef}},
   eprint = {2606.09221},
  year = 2026,
  month = Jun,
  eprintclass = {math.PR}
}
\end{document}